\documentclass{amsart}
\usepackage{amsmath, amssymb}

\newtheorem{theorem}{Theorem}
\newtheorem{lemma}[theorem]{Lemma}
\newtheorem{proposition}[theorem]{Proposition}
\newtheorem{corollary}[theorem]{Corollary}

\theoremstyle{definition}
\newtheorem{definition}[theorem]{Definition}
\newtheorem{example}[theorem]{Example}

\theoremstyle{remark}
\newtheorem{remark}[theorem]{Remark}

\numberwithin{equation}{section}
\numberwithin{theorem}{section}
\newcommand\lbb[1]{\label{#1}}
\def\<{\langle}
\def\>{\rangle}
\def\ov{\overline}
\def\d{\partial}
\def\st{\; | \;}                               
\newcommand\encirc[1]{{}_{(#1)}}
\newcommand{\ti}[1]{\widetilde{#1}}
\newcommand{\wti}{\widetilde}
\def\Cset{\mathbb{C}}       
\def\Zset{\mathbb{Z}}       
\newcommand{\CC}{\mathbb{C}}
\newcommand{\ZZ}{\mathbb{Z}}
\newcommand{\fg}{\mathfrak{g}}
\newcommand{\fa}{\mathfrak{a}}

\def\de{\delta}

\def\la{\lambda}

\def\g{{\mathfrak{g}}}      

\def\ss{{\mathfrak{s}}}
\def\rr{{\mathfrak{r}}}
\def\gl{{\mathfrak{gl}}}
\def\sl{{\mathfrak{sl}}}
\def\so{{\mathfrak{so}}}

\def\cso{{\mathfrak{cso}}}

\def\F{\mathcal{F}}

\DeclareMathOperator{\Res}{Res} 

\DeclareMathOperator{\ad}{ad}

\DeclareMathOperator{\sd}{\ltimes}
\DeclareMathOperator{\Div}{div}

\DeclareMathOperator{\Lie}{Lie}

\DeclareMathOperator{\cder}{Cder}
\DeclareMathOperator{\cinder}{Cinder}
\DeclareMathOperator{\inder}{Inder}

\DeclareMathOperator{\Tor}{Tor}
\DeclareMathOperator{\Rad}{Rad}

\DeclareMathOperator{\rk}{rk}
\DeclareMathOperator{\ann}{Ann}
\DeclareMathOperator{\End}{End}

\DeclareMathOperator{\cur}{Cur}

\DeclareMathOperator{\vir}{Vir}
\DeclareMathOperator{\der}{Der}
\DeclareMathOperator{\cend}{Cend}
\DeclareMathOperator{\chom}{Chom}
\DeclareMathOperator{\gc}{gc}
\DeclareMathOperator{\im}{im}

\begin{document}
                                                                                
\title{Structure theory of finite Lie conformal superalgebras}
                                                                                
\author[D.~Fattori]{Davide Fattori}
\address{Dipartimento di Matematica Pura e Applicata, Universit\'a di 
Padova,
via Belzoni 7, 35131 Padova, Italy}
\email{fattori@math.unipd.it}

\author[V.~G.~Kac]{Victor G.~Kac}
\address{Department of Mathematics, MIT, Cambridge MA 02139, USA}
\email{kac@math.mit.edu}

\author[A.~Retakh]{Alexander Retakh}
\address{Department of Mathematics, MIT, Cambridge MA 02139, USA}
\email{retakh@math.mit.edu}
\thanks{The second and the third authors were partially supported by the 
NSF}


                                                                                

\maketitle

\section*{Introduction}

Lie conformal superalgebras encode the singular part of the operator
product expansion of chiral fields in two-dimensional quantum field theory
\cite{K1}.  On the other hand, they are closely connected to the notion of
a formal distribution Lie superalgebra $(\fg,\mathcal{F})$, i.e. a Lie
superalgebra $\fg$ spanned by the coefficients of a family $\mathcal{F}$
of mutually local formal distributions. Namely, to a Lie conformal
superalgebra $R$ one can associate a formal distribution Lie superalgebra
$(\Lie R,R)$ which establishes an equivalence between the category of Lie
conformal superalgebras and the category of equivalence classes of formal
distribution Lie superalgebras obtained as quotients of $\Lie R$ by
irregular ideals \cite{K1}.

The classification of finite simple Lie conformal superalgebras was
completed in \cite{FK}.  The proof relies on the methods developed
in \cite{DK} for the classification of finite simple and semisimple Lie
conformal algebras, and the classification of simple linearly compact Lie
superalgebras \cite{K5}.

The main result of the present paper is the classification of finite
semisimple Lie conformal superalgebras (Theorem~\ref{mainth}).  Unlike in
\cite{DK}, we do not use the connection to formal distribution algebras in
the proof of this theorem (which would require us to take care of numerous
technical difficulties). We work instead entirely in the category of Lie
conformal superalgebras.  Our key result is the determination of finite
differentiably simple Lie conformal superalgebras 
(Theorem~\ref{diffsimple}).  The proof of this result uses heavily the 
ideas of \cite{B} and \cite{C}.

Given a finite Lie conformal superalgebra $R$, denote by $\Rad R$ the sum 
of all sol\-vable ideals of $R$.  Since the rank of a finite solvable Lie
conformal (super)algebra is greater than the rank of its derived
subalgebra \cite{DK}, we conclude that $\Rad R$ is the maximal solvable 
ideal, hence
$R/\Rad R$ is a finite semisimple Lie conformal superalgebra.  Thus in
some sense the study of general finite Lie conformal superalgebras reduces
to that of semisimple and solvable superalgebras.

In the Lie conformal algebra case there is a conformal analog of Lie's
theorem, stating that any non-trivial finite irreducible module over a
finite solvable Lie conformal algebra is free of rank $1$ \cite{DK}.  
However, a similar result in the ``super'' case is certainly false.  We
hope that we can develop a theory of finite irreducible modules over
solvable Lie conformal superalgebras similar to that in the Lie
superalgebra case \cite{K3}.  We make several observations to that end in
the last section of this paper.

The paper is organized as follows: we define the main objects of
our study and provide some general statements in
Section~\ref{sec.prelim}.  In Section~\ref{sec.diffsim} we
establish the structure of differentiably simple Lie conformal
superalgebras. Our proof is, in fact, quite general and the result
is valid for non-Lie finite conformal superalgebras as well.  We list
finite simple Lie conformal superalgebras in Section~\ref{sec.simple}
and describe their conformal derivations.  Then in
Section~\ref{sec.diffsimder} we describe conformal derivations in
the differentiably simple case and thus complete the
classification of finite differentiably simple Lie conformal
superalgebras.  This allows us to describe the structure of finite
semisimple Lie conformal superalgebras in
Section~\ref{sec.semisim}.  
In Section~\ref{sec.physpairs} we classify simple physical Virasoro
pairs and, as a consequence, obtain a
classification of physical Lie conformal superalgebras which generalizes 
that of \cite{K3}.  Finally,
in Section~\ref{sec.lieth} we initiate the study of
representations of finite solvable Lie conformal superalgebras.


\section{Basic definitions and structures}\label{sec.prelim}

\subsection{Formal distributions and conformal algebras}
Let $\g$ be a Lie superalgebra. A $\g$-valued \emph{formal
distribution} in one indeterminate $z$ is a formal power series
\begin{equation*}
a(z)=\sum_{n \in \Zset}a_{n}z^{-n-1}, \quad a_n \in \g.
\end{equation*}

The vector superspace of all formal distributions,
$\g[[z,z^{-1}]]$, has a natural structure of a
$\Cset[\d_z]$-module. We define
\begin{equation*}
\Res_z  a(z) = a_0.
\end{equation*}

Let $\g$ be a Lie superalgebra, and let $a(z),b(z)$ be two $\g$-valued
formal distributions. They are called \emph{local} if
\begin{equation*}
(z-w)^N [a(z),b(w)]=0 \quad \text{for} \quad N \gg 0.
\end{equation*}

Let $\g$ be a Lie superalgebra, and let $\F$ be a family of
$\g$-valued mutually local formal distributions. The pair $(\g,\F)$ is
called a \emph{formal distribution Lie superalgebra} if $\g$ is
spanned by the coefficients of all formal distributions from $\F$.

The bracket of two local formal distributions is given by the formula
\begin{equation*}
[a(z),b(w)]=\sum_j [a(w)\encirc{j}b(w)]\d_w^j\delta(z-w)/j!,
\end{equation*}
where $[a(w)\encirc{j}b(w)]=\Res_z(z-w)^j[a(z),b(w)]$.
Thus we get a family of operations
$\encirc{n}$, $n\in\ZZ_+$, on the space
of formal distributions: $[a(z)\encirc{n}b(z)]$ 
We define the \emph{$\lambda$-bracket} on the space
of formal distributions as the generating series of these operations 
\cite{DK,K1}:
\begin{equation*}
[a_{\lambda}b]=\sum_{n \in \Zset_+}\frac{\la^n}{n!}(a_{(n)}b).
\end{equation*}

The properties of the $\la$-bracket lead to the following basic
definition (see \cite{K1,K3}):

A \emph{Lie conformal superalgebra} $R$ is a left
$\Zset/2\Zset$-graded $\Cset[\d]$-module endowed with a $\CC$-linear
map, called the $\la$\emph{-bracket},
\begin{displaymath}
  R \otimes R \to \CC [\lambda]\otimes R ,\qquad a \otimes b \mapsto
  [a_{\lambda}b]
\end{displaymath}
satisfying the following axioms $(a,b,c \in R)$:

\begin{align*}
\begin{array}{ll}
    \text{(sesquilinearity)}\quad &
    [\partial a_{\lambda}b]=-\lambda [a_{\lambda}b], \;\;
    [a_{\lambda}\partial b]=(\partial +\lambda)[a_{\lambda}b],
    \\
    \text{(skew-commutativity)}\quad &
    [b_{\lambda}a]=-(-1)^{p(a)p(b)}[a_{-\lambda -\partial}b] ,\\
\text{(Jacobi identity)}\quad &
    [a_{\lambda} [b_{\mu}c]] =
    [[a_{\lambda}b]_{\lambda +\mu}c]+ (-1)^{p(a)p(b)}
   [b_{\mu}[a_{\lambda}c]] .
\end{array}
\end{align*}
(Here and further $p(a)\in\ZZ/2\ZZ$ stands for the parity of an element 
$a$.)

Similarly, one can define an \emph{associative conformal
superalgebra} by replacing the skew-commutativity and the Jacobi
identity above with the following property:
\begin{equation*}
\text{(associativity)}\qquad\quad a_{\lambda} (b_{\mu}c) =
    (a_{\lambda}b)_{\lambda +\mu}c.
\end{equation*}

As a shorthand, we call $R$ \emph{finite} if $R$ is finitely generated as 
a $\CC[\d]$-module.

Below follow two useful constructions of Lie conformal superalgebras 
(more concrete examples will be discussed in Section~\ref{sec.simple}):

\begin{example}\label{assoctolie} An associative conformal superalgebra
can be endowed with
a $\la$-bracket by putting
\begin{equation*}
[a_\la b]=a_\la b -(-1)^{p(a)p(b)}b_{-\lambda -\partial}a.
\end{equation*}
\end{example}

\begin{example}\lbb{tplcsos} Let $R$ be a Lie conformal superalgebra and
let $B$ be a commutative associative (ordinary) superalgebra.
Then $R\otimes B$ carries a Lie conformal superalgebra structure defined
as follows. The $\Cset[\d]$-module structure is given by
$\d(r\otimes b)= (\d r)\otimes b$ ($r\in R,\; b\in B$), and
the $\la$-bracket by
\begin{equation*}
[(r\otimes b)_{\la}(r' \otimes b')]=(-1)^{p(b)p(r')}
[r_{\la}r']\otimes(bb').
\end{equation*}
Notice that if $R$ is finite and $B$ is finite-dimensional, then $R
\otimes B$
is also finite.
\end{example}

\subsection{Structural terminology} We will denote by $\<X_\la Y\>$ (or
$\<[X_\la Y]\>$, whenever appropriate) the $\CC[\d]$-module generated by
elements of the form $x\encirc{n}y$, where $x\in X, y\in Y$, and 
$n\in\ZZ_+$.

An \emph{ideal} of a Lie conformal superalgebra $R$ is a 
$\CC[\d]$-submodule
$I$ of $R$ such that $\<[R_\la I]\>\subset I$.  An ideal $I$ is
\emph{abelian} if $[I_\la I]=0$, \emph{central} if $[R_\la I]=0$, and
\emph{nilpotent} if $\<[\ldots\<[I_{\la_1} I]\>_{\la_2}\ldots_{\la_k}
I]\>=0$
(in this case we sometimes write $I^{k+1}=0$).

A Lie conformal superalgebra is \emph{simple} if it is non-abelian and 
contains 
no ideals except for zero and itself.

The derived series of a Lie conformal superalgebra is built in the usual
fashion: let $R'=\<[R_\la R]\>$ and set $R^{(0)}=R$,
$R^{(n+1)}=(R^{(n)})'$.  Then a Lie conformal superalgebra $R$ is
\emph{solvable} if $R^{(n)}=0$ for some $n$. $\Rad R$ is the maximal
solvable ideal of $R$ (its existence is explained in the introduction).

A Lie conformal superalgebra is \emph{semisimple} if it has no non-zero
abelian ideals. Since the second last term of the derived series is
an abelian ideal, this is equivalent to saying that $R$ has no non-zero
solvable ideals. Since the $\Cset[\d]$-torsion is central \cite{DK},
a finite semisimple Lie conformal superalgebra is free
as a $\Cset[\d]$-module.

\subsection{Conformal modules} Let $V, W$ be  $\Zset_2$-graded left
$\Cset[\d]$-modules. We denote by
$\End_{\Cset[\d]} V$ the set of all $\Cset[\d]$-linear endomorphisms of
$V$.
Notice that $\End_{\Cset[\d]} V$
has a $\Cset[\d]$-module structure given by
\begin{equation*}
(\d \varphi)(v)= \d \varphi(v) \quad \text{for any} \quad v \in V.
\end{equation*}

A $\Cset$-linear map
\begin{equation*}
\varphi \; : \; V\to \Cset[\la] \otimes_{\Cset} W
\end{equation*}
is called a \emph{conformal linear map} if the following equation holds:
\begin{equation*}
\varphi_{\la}(\d v)=(\d +\la)\varphi_{\la}v \quad \text{for any }
v \in V.
\end{equation*}
The $\Cset$-vector space of all conformal linear maps from $V$ to $W$ is
denoted by $\chom(V,W)$.
It  has $\Cset[\d]$-module structure if we set
\begin{equation*}
(\d \varphi)_{\la}(v)=-\la \varphi_{\la} v.
\end{equation*}

When $V=W$, we denote $\chom(V,V)$ by $\cend V$. When $V$ is finite over
$\CC[\d]$, $\cend V$ becomes an associative conformal superalgebra with
the $\lambda$-product
\begin{equation*}
(a_\la b)_\mu v=a_\la(b_{\mu-\la}v),\quad a,b\in\cend V.
\end{equation*}
The space $\cend V$ endowed with a $\la$-bracket (see 
Remark~\ref{assoctolie}) is
denoted $\gc V$.

A \emph{module} $M$ over a Lie conformal superalgebra $R$ is a
$\Zset_2$-graded left $\Cset[\d]$-module $M$ endowed with a
$\Cset$-linear map
\begin{equation*}
R\to\gc M.
\end{equation*}

Alternatively, one can define a module over $R$ by providing a
map $R\to\CC[\la]\otimes\End_\CC M$, $a\mapsto a_\la^M$ such that
\begin{align*}
&(\d a)_{\la}^Mm=[\d,a_{\la}^M]m=- \la  a_{\la}^Mm,  \\
&[a^M_\lambda,b^M_\mu]v=[a_\la b]^M_{\la+\mu}v.
\end{align*}

An  $R$-module $M$ is \emph{simple} if it has no
nontrivial $R$-invariant $\Cset[\d]$-submodules.

An \emph{endomorphism} of an $R$-module $M$ is a
$\Cset[\d]$-linear map $\phi \in \End_{\Cset[\d]} M$ such that for any $a
\in R$ and $v \in M$
we have
\begin{equation*}
\phi(a_{\mu}v)=(-1)^{p(a)p(\phi)}a_{\mu}\phi(v).
\end{equation*}

\subsection{Derivations} A \emph{conformal derivation} of a Lie conformal
superalgebra $R$ is a conformal endomorphism $\phi$ of $R$ such that
for any homogeneous $x,y \in R$
\begin{equation*}
\phi_{\la}[x_{\mu}y]=[(\phi_{\la}x)_{\la + \mu}y] +
(-1)^{p(x)p(\phi)}[x_{\mu}
(\phi_{\la}y)].
\end{equation*}
We denote by $\cder R$ (resp. $\cinder R$) the space of all
conformal derivations of $R$ (resp. of all \emph{inner conformal derivations}, i.e.
conformal derivations
of the form $\ad a$, $a \in R$, $(\ad a)_\la b=[a_\la b]$ for $b\in R$).  
Clearly $\cder R$ is a subalgebra
of $\gc R$.

Let $D$ be a set of conformal derivations of $R$. An ideal $I$ of $R$ is
{\em $D$-stable} if $\phi I\subseteq I$ for all $\phi\in D$.  Here (and 
below) we use the shorthand $\phi I$ for $\<\phi_\la I\>$.

A Lie conformal superalgebra is {\em $D$-differentiably simple} if it 
contains no proper
$D$-stable ideals.  A Lie conformal superalgebra is {\em differentiably
simple} if it is $D$-differentiably simple with respect to some $D$.

We will also work with ordinary derivations of conformal algebras.  An
\emph{ordinary derivation} of a Lie conformal superalgebra $R$ is
a $\Cset[\d]$-linear endomorphism $d$ of $R$  such that for any 
homogeneous $x,y\in R$,
\begin{equation*}
d([x_{\mu}y])=[d(x)_{ \mu}y] + (-1)^{p(x)p(d)}[x_{\mu}
d(y)].
\end{equation*}
We denote the space of all ordinary derivations of $R$ by $\der R$.
Remark that a conformal derivation $\phi$ gives rise to an ordinary derivation
$\phi\encirc{0}$. In particular, every element $a\in R$ gives rise to an 
ordinary derivation $\ad a\encirc{0}$.  We call such derivations 
\emph{inner} 
and denote their space as $\inder R$.

\begin{remark} The operator $\d$ always acts as an ordinary derivation of 
a conformal algebra.  In some cases it is inner (e.g. when the algebra 
possesses a Virasoro element, see Example~\ref{vironcur}). 
\end{remark}

\subsection{Minimal ideals}  Here we collect some facts about minimal
ideals of finite Lie conformal superalgebras.

\begin{lemma}\lbb{milcs} Let $R$ be a finite Lie conformal superalgebra
and
$J$ an ideal of $R$ which contains no nonzero central elements.
Then $J$ contains a minimal  ideal of $R$.
\end{lemma}

\begin{proof} Let $I$ be a minimal rank ideal of $R$ contained in $J$.
Let $I_0=\cap_i ~K_i$ be the intersection of all non-zero
ideals of $R$ contained in $I$. Any ideal $K_i$ has the same rank as $I$,
hence
$I/K_i$ is a torsion $\Cset[\d]$-module and by \cite[Proposition 3.2]{DK},
$R$ acts trivially on it.
This means that for any $i$ we have
$\langle [R_{\la}I] \rangle \subseteq K_i$.
Therefore $\langle [R_{\la}I ]\rangle \subseteq I_0$.
Note that $\langle [R_{\la}I ] \rangle \neq 0$ because otherwise $I$
would be a central ideal of $R$. Hence $I_0 \neq 0$
and clearly $I_0$ is a minimal ideal of $R$.
\end{proof}

\begin{lemma}\lbb{fctsmi}
Let $M$ be a nonabelian ideal in a finite Lie conformal superalgebra $R$.
Then
\begin{enumerate}
\item
$M$ is a minimal ideal if and only if $M$ is $\ad R$-simple.
\item
If $M$ is a minimal ideal, then it is differentiably simple.
\item
If $M$ is a minimal ideal, then it is $\cder R$-invariant and $\cder 
R$-simple.
\item
If $M$ is  a minimal ideal, then either $M$ is simple or $M$
contains a minimal ideal $I$ of $M$ which is abelian.
\end{enumerate}
\end{lemma}

\begin{proof}
(1) is immediate. (2) follows from (1) and the fact that $\ad R
\subseteq
\cder M$. As for (3), we remark that the minimality of $M$ implies that
$\langle [M_{\la}M] \rangle =M$, hence $ \langle \phi_{\la}(M) \rangle
\subseteq \<[(\phi_\mu M)_\la M]\> \subseteq M$
for any $\phi \in \cder R$, i.e. $M$ is a $\cder R$-invariant ideal of
$R$.
Now, let $J$  be a nonzero $\cder R$-invariant ideal
of $R$, which is contained in $M$. The minimality of $M$ implies that
$J=M$.
In order to prove (4), notice that by
(2) $M$ is differentiably simple. Suppose $M$ is not simple.
The center of $M$ is a differential ideal of $M$, hence it is zero and
Lemma \ref{milcs} provides a minimal ideal $I$ in $M$. Suppose that $I$ is
not abelian. Then $\langle [I_{\la}I] \rangle =I$ and $I$ is
a $\cder M$-invariant ideal in $M$ which is differentiably simple. Hence
$I=M$ i.e. $M$ is a minimal ideal in itself. We conclude that
$M$ is simple. The contradiction proves that $I$ is abelian.
\end{proof}

Note that by Lemma~\ref{fctsmi}(3), any minimal ideal in a differentiably
simple but non-simple finite Lie conformal superalgebra is abelian.


\section{Differentiably simple conformal
superalgebras}\label{sec.diffsim}

In this section we prove the following

\begin{theorem}\label{diffsimple} Let $R$ be a finite 
non-abelian differentiably simple Lie conformal superalgebra.  Then 
$R\simeq S\otimes\wedge(n)$ for a simple Lie conformal superalgebra $S$.
\end{theorem}

\begin{remark}Our proof also works for any finite 
non-abelian differentiably simple conformal
superalgebra (i.e. the one for which only sesquilinearity holds) but we do 
not require the result in this generality.  
\end{remark}

\subsection{Centroid} Let $R$ be a Lie conformal superalgebra and $M$
an $R$-module. For $x\in R$, let $L_x$ be an element of $\cend R$
such that $({L_x})_\lambda y=x_\lambda y$ for any $y\in M$.

\begin{definition} The {\em centroid} $C(M)$ of a module $M$ over a Lie
conformal superalgebra $R$ is the subalgebra of the 
associative superalgebra $\End_\CC
M$ that consists of elements (super)commuting with $(L_x)\encirc{n}$
for all $x\in R$, $n\in\ZZ_+$ and the action of $\d$.
\end{definition}

\begin{remark}  By definition $C(M)$ is a subalgebra of the 
associative superalgebra
$\End_{\CC[\d]} M$.  
\end{remark}

\begin{remark} We We show in the proof of Lemma~\ref{commcent} below that 
for $a\in C(R)\subset 
\End_{\CC[\d]} R$, $a(x_\la y)=(ax)_\la y=(-1)^{p(a)p(x)}(x_\la ay)$.  
These conditions can be taken as the definition of the centroid of a 
conformal superalgebra $R$.  (Note that here $R$ is not necessarily 
Lie.)
\end{remark}

\begin{lemma}\label{commcent} If $R=\<[R_\la R]\>$, then $C(R)$ is 
(super)commutative.
\end{lemma}

\begin{proof}
The lemma follows from the equalities $a[x_\lambda
y]=(-1)^{p(a)p(x)}[x_\lambda ay]$ and $a[x_\lambda y]=[ax_\lambda
y]$ for any $a\in C(R), x,y\in R$.  The first equality follows directly 
from the definition; the second is deduced from the first:
\begin{align*}
a[x_\lambda y]=&-(-1)^{p(x)p(y)}a[y_{-\la-\d} x]= 
-(-1)^{p(x)p(y)+p(a)p(y)}[y_{-\la-\d} ax]\\
=&(-1)^{p(x)p(y)+p(a)p(y)+(p(a)+p(x))p(y)}[ax_\la y]=[ax_\la y].
\end{align*}
Then for any $a,b\in C(R)$ and any $x,y\in R$, $(ab)[x_\la 
y]=(-1)^{p(b)p(x)}[ax_\la by]$ and $(ba)[x_\la 
y]=(-1)^{(p(a)+p(x))p(b)}[ax_\la by]$, implying that $a$ and $b$ 
(super)commute. 
\end{proof}

Since $\<[R_\la R]\>$ is a differentiably stable ideal, it follows
that $C(R)$ is (super)com\-mutative for a differentiably simple
$R$.

\begin{lemma} For a homogeneous (i.e. even or odd) $a\in C(R)$, $\ker a$ 
and $\im a$ are ideals of $R$.
\end{lemma}

\begin{proof} Clear.
\end{proof}

We also have a version of the Schur Lemma:

\begin{lemma}\lbb{cvsclmlcs}
Let $M$ be a countable dimensional simple module over a Lie conformal
superalgebra $R$. Then either $C(M)=\Cset 1_M$ or $\dim M_{\ov{0}}=
\dim M_{\ov{1}}$ and
\begin{equation*}
C(M)=\Cset 1_M  \oplus \Cset  U,
\end{equation*}
where $U$ is an odd operator such that $U^2=1_M$.
\end{lemma}

\begin{proof} The proof follows the classical line of argument.
Let $a \in C(M)$ be a non-zero even operator. The fact that $M$ is
simple implies that $a$ is invertible. Suppose $a$ is not a
scalar. Since $\Cset$  is algebraically closed, $a$ cannot be an
algebraic element in $C(M)$. The field of rational functions 
in $a$ over $\CC$ is contained in $C(M)$, hence
$C(M)$ has dimension greater than countable.
On the other hand, let us fix a nonzero $x \in M$. The map $C(M)
\to M$ sending $a$ to $ax$ is injective because $C(M)$ is a
division ring. However, $M$ is countable dimensional. The
contradiction proves that $a =c 1_M$ for some $c \in \Cset$.

Let $a\in C(R)$ be a non-zero odd operator. Then $b$ is invertible
(this can only happen if $\dim M_{\ov{0}}= \dim M_{\ov{1}}$).
Furthermore, $a^2$ is an even operator, hence a scalar. Suppose
$a_1, a_2$ are two non-zero odd operators in $C(M)$. Then $(a_1-c
a_2)^2\in\CC 1_M$ for any $c\in \Cset$.  Assume that for $c$ and $c'$,
$(a_1-c a_2)^2\neq 0$ and $(a_1-c'a_2)^2\neq 0$.  Then $(a_1-c a_2)^2
+r^2(a_1-c'a_2)^2=0$ for some $0\neq r\in\CC$.  We obtain two non-zero
odd operators $b_1$ and $b_2$ such that $b_1^2+b_2^2=0$.  By taking the
square of $b_1b_2b_1^{-1}b_2^{-1}$ one easily shows that $b_1b_2=\pm
b_1b_2$.  Thus either $(b_1+b_2)^2=0$ or $(b_1+ib_2)(b_1-ib_2)=0$ and
it follows that $a_1$ is proportional to $a_2$.
\end{proof}

\begin{corollary}\label{schursim} If $R$ is a countably dimensional simple 
Lie conformal superalgebra, then $C(R)=\CC 1_R$.
\end{corollary}

\begin{proof}  Assume $C(R)=\CC 1_R\oplus \CC U$, $U^2$=1.  Let
$a=(1+U)/2$ and $\bar a=(1-U)/2$.  $R$ splits as $aR\oplus \bar a R$ and
it is easy to see that $\ker a=\bar a R$.
Using the equality $b[x_\la
y]=[bx_\la y]$, we see
immediately that $\ker a$ and $\im a$ are (non-homogeneous) ideals of $R$.
On the other hand, $[x_\la ay]=\bar a[x_\la y]$ for $x$ odd, hence $aR\cap
\bar a R\neq\{0\}$, a contradiction.
\end{proof}

\begin{remark}\label{schurlie} Lemma~\ref{cvsclmlcs} and 
Corollary~\ref{schursim} hold for ordinary Lie superalgebras (with the 
same proof).
\end{remark}

\subsection{Conformal Derivations and the centroid}   Let $\phi\in\cder 
R$, $\gamma\in\CC$, and $n\in\ZZ_+$. We introduce the 
following operators acting on $\End R$: 
$\ad\phi_\gamma=[\phi_\gamma,\bullet]$ and 
$\ad\phi\encirc{n}=[\phi\encirc{n},\bullet]$.

\begin{lemma} Let $a\in C(R)$, $\phi\in\cder R$.

1) For any $\gamma\in\CC$, $\ad\phi_\gamma (a)\in C(R)$.  Moreover,
$\ad\phi_\gamma$ is a derivation of $C(R)$.

2) For any $n\in\ZZ_+$, $\ad\phi\encirc{n}(a)\in C(R)$.
Moreover, $\ad\phi\encirc{n}$ is a derivation of $C(R)$.
\end{lemma}

\begin{proof} Direct computation.
\end{proof}

\begin{lemma} If $a\in C(R)$ and $\phi\in Cder R$, then
$a\phi\in Cder R$. Also $p(a\phi)=p(a)+p(\phi)$.
\end{lemma}

\begin{proof} Direct computation.
\end{proof}

\subsection{Constructing a chain of ideals} Let $R$ be a finite 
non-abelian
differentiably simple Lie conformal superalgebra.  Remark that the center
of $R$, $\{x\st [x_\la R]=0\}$ is a differentiably stable ideal, hence 
it is zero and by
Lemma~\ref{milcs}, $R$ contains a minimal ideal $I$.  Let $D$ be a
set of homogeneous conformal derivations of $R$. Our ultimate goal is to
construct with the use of $D$ a certain finite chain of ideals that starts 
with $I$ (cf. \cite{B}).

Assume now that we have constructed a chain of ideals
$I=I_1\subset I_2\subset\dots\subset I_q\neq R$ such that for all
$2\leq j\leq q$, $I_j/I_{j-1}\simeq I$ as $R$-modules.  Let
$\phi\in D$ be a homogeneous conformal derivation such that $\phi
I_q\not\subset I_q$. We are going to construct an ideal
$I_{q+1}\supset I_q$ such that $I_{q+1}/I_q\simeq I$ as
$R$-modules.

Remark, first that for any ideal $J$, a homogeneous conformal derivation $\phi$
induces
a map $\bar\phi\in\chom(J, R/J)$, $\bar\phi_\lambda x=\phi_\lambda
x+\CC[\lambda]\otimes J$.  Moreover, by definition, for every
$y\in R$, the following equalities holds in the $R$-module $R/J$:
\begin{gather*}
\bar\phi_\lambda(y_\mu x)=(-1)^{p(\phi)p(y)}(y_\mu (\bar\phi_\lambda
x));\\
\d(\bar\phi_\lambda x)=\bar\phi_\lambda(\d
x)-\lambda\bar\phi_\lambda x.
\end{gather*}
It follows that for any $\gamma\in\CC$, $\ker\bar\phi_\gamma$ is a
homogeneous ideal of $R$ and $\im\bar\phi_\gamma$ is an
$R$-submodule of $R/J$.

Suppose now that we have constructed a chain such as above:
$I_1\subset\dots\subset I_q\neq R$.  Let $j$ be minimal such that
$\phi I_j\not\subset I_q$. Take $J=I_q$ and consider the map
$\bar\phi\in\chom(I_q,R/I_q)$ constructed as above. Restrict this
map to $I_j$.  By construction $\phi(I_{j-1})\subseteq I_q$, so
$I_{j-1}\subseteq \ker\bar\phi_\gamma$ for any $\gamma\in\CC$.  Thus
we have a family of maps $\bar\phi_\gamma:I_j/I_{j-1}\to R/I_q$.

We will show that there exists $\gamma$ such that
$\ker\bar\phi_\gamma$ is zero.  Indeed, if for some $\gamma\in\CC$
there exists $x+I_{j-1}\in I_j/I_{j-1}$ such that $\bar\phi_\gamma
x=0$, then $\bar\phi_\gamma (I_j/I_{j-1})=0$ (it is a simple
$R$-module by construction).  Hence, if there exists such an
$x+I_{j-1}$ for every $\gamma$, then
$\bar\phi_\lambda(I_j/I_{j-1})=0$ and, consequently, $\phi
I_j\subset I_q$.

Hence we can find $\gamma\in\CC$ with $\ker\bar\phi_\gamma=0$.
Then we can define $I_{q+1}$ as an ideal such that
$I_{q+1}/I_q=\im\bar\phi_\gamma$.

\begin{remark} If we start constructing a chain of ideals such as above
and keep the same $\phi\in D$, we would at some point obtain $I_q$
such that $\phi I\subset I_q$.
Indeed, pick a basis $x_1,\dots,x_n$ of $I$ over $\CC[\d]$ and let
$d$ be the maximal degree in $\lambda$ of all $\phi_\lambda x_i$.
At every step of the above construction, we produce an ideal $I_k$
and $\gamma_k\in\CC$ such that $\phi_{\gamma_k} I_1\subset I_k$.
It is clear that $\gamma_1,\dots,\gamma_d$ are pairwise distinct,
thus $\phi_\lambda I_1\subset I_d$ as required.
\end{remark}

\subsection{Constructing the maximal ideal}  Now let $D$ be a finite
collection of homogeneous conformal derivations,
$D=\{\phi_1,\dots,\phi_m\}$, such that $I$ is not $D$-stable.  Let
$i_1$ be the least index for which $\phi_{i_1}I\not\subset I$.
We apply the algorithm from the previous subsection.
Namely, using $\phi_{i_1}$ we can construct a chain of ideals
$I=I_1\subset\dots\subset I_{r_1}$, where $\phi_{i_1}I\subset
I_{r_1}$. Then use $\phi_{i_2}$ (with the minimal $i_2$) such that
$\phi_{i_2}I_{r_1}\not\subset I_{r_1}$ to extend the chain to
$I_{r_2}$, where $\phi_{i_2}I_{r_1}\subset I_{r_2}$, etc.

Either at some point we obtain $I_l=R$ or we obtain a proper ideal
$I_q$ such that $\phi_m\dots\phi_1 I\subseteq I_q$.

Suppose now that for every collection $D$ only the second case
occurs. Since $R$ is differentiably simple, there exists a homogeneous 
conformal derivation $\phi$ such
that $\phi I_q\not\subset I_q$ and, using $\phi$, we can extend
our chain to $I_{q+1}$.  For any $x\in I$, let
$y=\phi_m\encirc{n_m}(\dots\phi_1\encirc{n_1}x)\dots)$.  Hence
$y\in I_q$ and $y_\lambda I_1\simeq y_\lambda (I_{q+1}/I_q)=0$.
But every $y\in R$ can be expressed as a $\CC[\d]$-linear
combination of elements in the above form (for some collection
$D=\{\phi_1,\dots,\phi_m\}$).  Hence, $[R_\lambda I]=0$, i.e. $R$ has a 
non-trivial center, a contradiction.

Therefore, there exists a finite chain $I=I_1\subset\dots\subset
I_l=R$ such that for every $2\leq j\leq l$, $I_j/I_{j-1}\simeq I$
as $R$-modules. Denote $I_{l-1}$ as $N$.  Then
$R/N\simeq I$ and $N$ is maximal.

Since $N_\lambda (R/N)=0$, $\<[N_\lambda I_j]\>\subset I_{j-1}$ for
any $j$. Thus $N^l=0$ and $N$ is nilpotent.

Suppose that there exists another maximal ideal $N'$ of $R$.  Then
$N+N'=R$ and $R/N'$ is nilpotent and simple, a contradiction (as
$R=\<[R_\la R]\>$). We arrive at

\begin{proposition}\label{maxid}  Let $R$ be a finite 
non-abelian differentiably
simple conformal Lie superalgebra and let $I$ be its minimal ideal.  
Then
$R$ possesses a chain of ideals $I=I_1\subset
I_2\subset\dots\subset N\subset R$, where $N$ is a unique maximal
ideal of $R$.  Each factor $I_j/I_{j-1}$ is isomorphic to $I$ as
$R$-modules. Moreover, $N$ is nilpotent and $R/N\simeq I$.
\end{proposition}

Actually, above we only used that $I$ is finite, this the proof implies 
the following aside corollary.

\begin{corollary}  If $R$ is a non-abelian differentiably simple Lie 
conformal superalgebra with a finite minimal ideal, then $R$ itself is 
finite.
\end{corollary}

\begin{remark}  Uniqueness of $N$ implies that all minimal ideals
of $R$ are $R$-isomorphic.
\end{remark}

\subsection{Centroid structure}  Here we construct an embedding $\sigma:
C(R/N)\to C(R)$. For $a\in C(R/N)$, define $\sigma(a)$ as the
composition of maps
\begin{equation*}
R\twoheadrightarrow R/N \stackrel{a}{\to} R/N \stackrel{\rho}{\to}
I\hookrightarrow R
\end{equation*}
(here the isomorphism $\rho:R/N \stackrel{\sim}{\to} I$ is the one
constructed in the proof of Proposition~\ref{maxid}).

Clearly, $\sigma$ is an embedding; moreover, $\im\sigma(a)=I$ for
all $a$.

For an arbitrary $R$-isomorphism
$\theta:R/N\stackrel{\sim}{\to}I$, consider the map
$\beta:R/N\stackrel{\theta}{\to}I\stackrel{\rho^{-1}}{\to} R/N$.
Since $\beta$ commutes with the action of $R$ and $\d$, we see
that $\sigma(\beta)=\theta$.

Let $I=I_1\subset I_2\subset\ldots\subset N\subset R$ be the chain
described in Proposition~\ref{maxid}.

\begin{lemma}\label{sigmaprop} There exists a family of monomorphisms 
$\sigma_q: C(R/N)\to
C(R)$, $1\leq q\leq l$, such that for every $a\in C(R/N)$
$\sigma_q(a)R+I_{q-1}=I_q$ and $\sigma_q(a)N\subseteq I_{q-1}$.
Moreover, if $\theta: R/N\to I_q/I_{q-1}$ is an $R$-isomorphism,
then there exists $b\in C(R/N)$ such that $\theta$ is induced by
$\sigma_q(b)$.
\end{lemma}

\begin{proof}  We put $\sigma_1=\sigma$ constructed above.  Then denote by
$d_q$ the map $\phi_\gamma$ that we used to construct $I_{q+1}$
and let $I_j$ be the ideal used in that construction.  Put
$\sigma_{q+1}=[d_q,\sigma_j]$ (a super-bracket).  Modulo
$I_{j-1}$, $\sigma_j(a)R=I_j$ for any nonzero $a$, hence modulo
$I_q$, $d_q\sigma_j(a)R=I_{q+1}$.  Since $\sigma_j(a)d_qR\subseteq
I_j\subset I_q$, $\sigma_{q+1}R+I_q=I_{q+1}$ as required.
Similarly, $\sigma_{q+1}N\subseteq I_q$.

Given $\theta:R/N\to I_{q+1}/I_q$, we can extend it to a map
$R/N\to I_j/I_{j-1}$ that is induced by $\sigma_j(b)$ for some
$b$.  Thus $\sigma_{q+1}(b)$ induces $\theta$.
\end{proof}

Using the maps $\sigma_q$ constructed in the above lemma, we have

\begin{lemma}\label{sigmamap} The map $\ti\sigma=\oplus_{q=1}^l\sigma_q$ 
is an
isomorphism $\displaystyle{\bigoplus_{\text{$l$ copies}}}
C(R/N)\stackrel{\sim}{\to} C(R)$.
Also, $J=\sigma_1(C(R/N))$ is a
minimal ideal of $C(R)$.
\end{lemma}

\begin{proof} Let $a\in C(R)$ and $q$ be minimal such that $aR\subseteq
I_q$.  Combined with a projection $R\to R/I_{q-1}$, $a$ induces a
map $R\to I_q/I_{q-1}$.  Its kernel is a maximal ideal, i.e. $N$.
Hence $a$ induces an  isomorphism $R/N\to I_q/I_{q-1}$.  Then
there exists $b\in C(R/N)$ such  that $(a-\sigma_q(b))R\subseteq
I_{q-1}$.  By induction, $\ti\sigma$ is onto.

Conversely, if $\ti\sigma(\oplus_1^l a_i)=0$, then
$\sigma_l(a_l)R\subseteq N$.  Thus $a_l=0$.  By induction, all
$a_i=0$.

The set $J=\sigma_1(C(R/N))=\{a\,|\,aR\subset I\}$ is obviously an
ideal of $C(R)$.  Let $a',a''\in J$.  Then they induce
isomorphisms $\theta',\theta'': R/N\to I$.  Then
$a'=ba''$, where $b\in C$ induces $\theta'^{-1}\theta''$.
\end{proof}

\begin{remark} We actually obtain a chain of ideals
$J=J_1\subset\dots\subset J$ with properties similar to that of
$I_1\subset\dots\subset R$.
\end{remark}

\begin{corollary} $C(R)$ is finite-dimensional.
\end{corollary}

\begin{proof} Follows from Lemmas~\ref{cvsclmlcs} and~\ref{sigmamap}.
\end{proof}

Denote by $D$ the subset of $\der C(R)$ that consists of
derivations $\ad\phi_\gamma$, $\phi\in\cder R$, $\gamma\in\CC$.

\begin{lemma}\label{diffsimcent} $C(R)$ is $D$-differentiably simple.
\end{lemma}

\begin{proof} Let $H$ be a $D$-stable ideal of $C(R)$.  Then
$\phi_\gamma HR\subseteq HR$ for any $\phi\in\cder R$,
$\gamma\in\CC$. Hence, $\phi(HR)\subseteq HR$ and $HR=R$.  There
exists $h\in H$ such that $hR+N=R$ (as $N$ is maximal).  Thus
$0\neq J_1h\subseteq H\cap J_1$, i.e. $J_1\subseteq H$. The construction 
of maps $\sigma_q$ and $\ti\sigma$ in Lemmas~\ref{sigmaprop} 
and~\ref{sigmamap} implies that $C\subseteq H$.
\end{proof}

\begin{corollary}\label{diffsimcentcor} Let $D'=\{\ad \phi\encirc{n}\}$.
Then $C(R)$ is $D'$-differentiably simple.
\end{corollary}

\begin{proof} A $D'$-stable ideal of $C(R)$ is $D$-stable
by the definition of $\phi_\gamma$.
\end{proof}

Since $C(R)$ is differentiably simple and finite-dimensional, we
obtain from \cite{C} the following

\begin{proposition} Let $C(R)$ be the centroid of a finite 
non-abelian differentiably simple conformal Lie superalgebra.  Then $C(R)$ 
is a Grassmann superalgebra $\wedge(r)$.
\end{proposition}

\subsection{More on differential simplicity} At this point the
non-conformal argument for a Lie superalgebra proceeds as follows: first, 
establish that
$C(L)$ is $\ad d$-differentiably simple for just one conformal derivation 
$d\in\der L$
(with certain additional conditions on $d$ in the super case)  and
then conclude that $L$ is $d$-differentiably simple.

In order for this to work for a Lie conformal superalgebra $R$ as well, we 
need conformal analog of $d$ of the form $\ad\phi_\gamma$ for 
some $\gamma$.  We need to tweak the proof of Theorem~5.1 in
\cite{C}.

We restate the theorem itself first.  Let $C$ be a Grassmann
algebra over an algebraically closed field of characteristic $0$.
Let $C$ be $D$-differentiably simple with respect to a homogeneous set of 
derivations $D$
which is both a subalgebra and a left $C$-module.  Then there
exists $d\in D$ such that $C$ is $d$-differentiably simple and for the
corresponding chain $J=J_1\subset\dots\subset C$,  the map $\bar
d: J_q\to C/J_q$ is homogeneous (i.e. $dx\equiv d_\epsilon x \mod
J_q$ for all $x\in J_q$, $\epsilon=\bar 0,\bar 1$).

The reason for the last condition is simple: in order to build a
chain of ideals we need to use a homogeneous derivation at every
step.  In this case we say that $d$ is {\em homogeneous at every
step}.

The major steps of the proof are: find a homogeneous nilpotent $m\in
C=\wedge(r)$ and a derivation $d_1$ such that $d_1(m)=1$.  Then
$D=D_0\oplus Cd_1$, where $D_0=\{d-d_1(m)d_1\,|\,d\in D\}$ (in
particular, $D_0$ acts as $0$ on $m$).  View $C$ as a $\langle
C,D_0\rangle$-module and use $d_1$ to construct the appropriate
chain of ideals.  The minimal ideal in this chain is isomorphic to
$\wedge(r_1)$, $r_1<r$, so we can use $d_0\in D_0$ to refine
it into the chain of ideals of $C$.

It follows from the proof that the requirement of $D$ being an
algebra is superfluous, thus in order to apply the construction of
$d_0$ and $d_1$ to the set $D'$ constructed in
Corollary~\ref{diffsimcentcor}, we have to show that it is a
$C$-module.  It suffices to show that
$d=\ad\phi\encirc{n}+\ad\psi\encirc{m}\in D'$ for $m\neq n$.  Let
$m<n$. Then
$d=(\ad\phi+(-1)^{n-m}\frac{m!}{n!}\d^{n-m}\psi)\encirc{n}$.

In the same vein we can assume that $d_0$ and $d_1$ are of the
form $\ad\phi\encirc{n}$ for the same $n$.  Hence, we can
strengthen the theorem as

\begin{lemma} There exist $\phi\in\cder R$ and $n$ such that
$C(R)=\wedge(r)$ is $\ad\phi\encirc{n}$-differentiably simple and
$\ad\phi\encirc{n}$ is homogeneous at every step.
\end{lemma}

\begin{lemma} There exist $\phi\in\cder R$ and $\gamma\in\CC$ such that
$C(R)$ is $\ad\phi_\gamma$-differentiably simple.
\end{lemma}

\begin{proof} Let $\phi$ and $n$ be as in the previous lemma.  Assume that
for each $\gamma\in\CC$, there exists an ideal $J_\gamma$ that is
$\ad\phi_\gamma$-stable.  Since $C(R)\simeq\wedge(n)$ has a
unique minimal ideal, $\cap J_\gamma$ is non-empty.  Thus, there
exists a proper ideal $J$ that is $\ad\phi_\gamma$-stable for all
$\gamma$.  Hence, $\ad\phi_\lambda$ maps $J$ into $J[\lambda]$ and
$J$ is $\ad\phi\encirc{n}$-stable.
\end{proof}

Clearly, $\ad\phi_\gamma$ is not necessarily homogeneous at every
step. Let $D$ be the $C(R)$-module generated by the homogeneous
components of $\ad\phi_\gamma$.  Notice that since
$\ad\chi_\gamma+\ad\phi_\gamma=\ad(\chi+\psi)_\gamma$, every
element of $D$ arises from a conformal derivation of $R$.

Since $C(R)$ is $D$-differentiably simple, we can apply the previous 
argument to obtain

\begin{lemma}\label{oneder-c} There exist $\psi$ and $\gamma$ such that
$C(R)$ is $\ad\psi_\gamma$-differentiably simple.  Moreover, 
$\ad\psi_\gamma$ 
is
homogeneous at every step.
\end{lemma}

Let $\psi$ and $\gamma$ be as in Lemma~\ref{oneder-c}.  We claim that $R$
contains no $\psi$-stable homogeneous ideals.  Indeed, let $M$ be
such an ideal.  We can always assume that $I\subset M$ (e.g. let
$I$ be a minimal ideal contained in $M$).  Moreover, $M\subset N$.
Let $H=\{a\,|\,aR\subset M\}$, a proper homogeneous ideal of
$C(R)$.  A direct calculation shows that $H$ is
$\ad\psi_\gamma$-stable, a contradiction.

Thus $R$ is $\psi$-stable.

We need to show that $\psi_\gamma$ acts homogeneously at every
step, i.e. that the map $\bar\psi_\gamma: I_q\to R/I_q$ is
homogeneous.  This will allow us to build the chain $\{I_q\}$. It
suffices to show that one of the homogeneous components of
$\bar\psi_\gamma$ acts by zero.  So, let $\phi$ and $\beta$ be
such that $\phi_\beta$ is homogeneous and $\ad\phi_\beta:J_q\to
C/J_q$ acts by zero.  (We assume here that we have already built
the chain $\{J_q\}$ and that the chain of ideals of $R$ has been
constructed up to $I_q$.  We also assume that $\phi_\beta I_{q-1}\subset
I_q$.) There exist $c\in J_q$ and $y\in R$ such that $cy$
generates $I_q/I_{q-1}$ over $R$.  Then $\psi_\beta
cy=[\psi_\beta,c]y\pm c\psi_\beta y\in I_q$. It follows that
$\psi_\beta R_\lambda cy\subset I_q[\lambda]$.  Summing up, we
have

\begin{lemma}\label{hstep} Let $R$ be a finite non-abelian differentiably
simple Lie
conformal algebra.
Then there exist a conformal derivation $\psi$ and 
$\gamma\in\CC$ such that $R$
is $\psi_\gamma$-differentiably simple.  Moreover, using $\psi_\gamma$ we 
can construct
a chain of homogeneous ideals $I=I_1\subset\dots\subset N\subset
R$.
\end{lemma}

\subsection{Splitting $R$}   Let $R$ be a finite non-abelian
differentiably simple Lie conformal algebra. Let $\psi$ and $\gamma$ be as
in the Lemma~\ref{hstep}.
Let $S=\{x\,|\,\psi_\gamma x\in I\}$.

\begin{lemma} $R=S\oplus N$ as $\CC[\d]$-submodules and $S$ is 
a conformal subalgebra of $R$.
\end{lemma}

\begin{proof} A direct computation shows that $S$ is closed with respect
to the $\d$-action and the $\lambda$-bracket inherited from $R$.

Let $0\neq x\in S\cap N$.  Then $x\in I_q\backslash I_{q-1}$ for
some $q$. Hence, $\psi_\gamma x\not\in I_q$, a contradiction, and
$S\cap N=0$.

To show that $R=S+N$, observe that $I_q\subseteq I+\psi_\gamma N$
for all $q$.  Indeed, by induction $I_{q-1}\subseteq I+\psi_\gamma
N$ and $\psi_\gamma I_{q-1}\subset \psi_\gamma N$.  Hence
$I_q=I_{q-1}+\psi_\gamma I_{q-1}\subseteq I+\psi_\gamma N$.  In
particular $R\subseteq I+\psi_\gamma N$.  Thus for every $x\in R$,
$\psi_\gamma x=\psi_\gamma y \mod I$ for some $y\in N$ and $x-y\in
S$.
\end{proof}

\begin{corollary} $S$ is a finite simple Lie conformal superalgebra.
\end{corollary}

\begin{proof} $S\simeq R/N$ as conformal superalgebras.
\end{proof}

Let $f:S\to R$ be the natural embedding. Define the map $F:
\wedge(r)\otimes S\to R$ by setting
\begin{equation*}
F(c\otimes x)=(-1)^{p(c)p(x)} cf(x),\qquad c\in \wedge(r),x\in R
\end{equation*}
for homogeneous elements.

Recall that we have constructed the chain $\{J_q\}$ of ideals of
$C(R)$. In particular, $J_q=\sigma_1(C(S))+\dots+\sigma_q(C(S))$.
Denote $c_q=\sigma_q(1)$ (it is homogeneous by definition of
$\sigma_q$). Clearly $F(c_1\otimes S)=I$.  By induction
$I_q+F(c_q\otimes S)=I_{q+1}$, hence $F$ is surjective.

To show injectivity, assume that $x$ is such that $F(x)=0$ and
$x=\sum_1^q c_i\otimes x_i$, $x_q\neq 0$.  Since $f(x_q)\not\in N$
and $c_q$ induces an isomorphism between $R/N$ and $I_q/I_{q-1}$,
$c_qf(x_q)\not\in I_{q-1}$.  On the other hand, $F$ maps the first
$q-1$ components into $I_{q-1}$.  Thus $c_qf(x_q)=0$.  Since
$c_qf$ is a homogeneous map, we obtain $x_q=0$.

This completes the proof of Theorem~\ref{diffsimple}.


\section{Simple Lie conformal superalgebras and their
derivations}\label{sec.simple} Here we recall the construction of all
simple Lie conformal superalgebras and describe their conformal 
and ordinary derivations.

\subsection{Simple conformal superalgebras}  In this subsection we review
the main result of \cite{FK}.

\begin{example}\lbb{loop}
Let $\g$ be a finite-dimensional Lie superalgebra. The \emph{loop
algebra} associated to $\g$ is the Lie superalgebra
\begin{equation*}
\wti{\g}= \g[t,t^{-1}], \;\;
p(at^k)=p(a) \;\; \text{for} \;\; a\in \g,\; k\in \Zset,
\end{equation*}
with the bracket
\begin{equation*}
[a\otimes t^n,b\otimes t^m]=[a,b]\otimes t^{n+m} \quad (a,b \in \g; \;m,n
\in
\Zset).
\end{equation*}
We introduce the family $\F_{\g}$ of formal distributions (known as
currents)
\begin{equation*}
a(z)=\sum_{n \in \Zset} (a\otimes t^n) \;z^{-n-1}, \quad a\in \g.
\end{equation*}
It is easily verified that
\begin{equation*}
[a(z),b(w)]= [a,b](w)\de(z-w),
\end{equation*}
hence $(\wti{\g}, \F_{\g})$ is a formal distribution Lie superalgebra.
The associated Lie conformal superalgebra is $\Cset[\d] \otimes \g$, with
the $\la$-bracket (we identify $1 \otimes \g$ with $\g$)
\begin{equation*}
[a_{\la}b]=[a,b], \quad a,b \in \g.
\end{equation*}
It is called the \emph{current conformal algebra} associated to $\g$,
and is denoted $\cur \g$.  The Lie conformal superalgebra
$\cur \g$ is simple if and only if $\g$ is a simple Lie superalgebra.
\end{example}

\begin{example}\label{vironcur} We define a conformal linear map $L: \cur
\g \to \cur \g$ by $L_{\la}g=(\d + \la)g $ for any $g \in \g$. It
is immediate to verify that this is a conformal derivation of
$\cur \g$. $L$ generates the conformal algebra $\CC[\d]L\subset
\gc(\cur \g)$ called the \emph{Virasoro conformal algebra} and
denoted $\vir$.  The $\lambda$-bracket is
\begin{equation*}
[L_\la L]=(\d+2\la)L.
\end{equation*}
Thus we have constructed a Lie conformal algebra $\vir\ltimes \cur\g$.

$\vir$ can be constructed using formal distributions with coefficients 
in the Lie algebra $\CC[t,t^{-1}]\d_t$, letting 
\begin{equation*}
L(z)=\sum_{n\in\ZZ} (t^n\d_t)z^{-n-1}.
\end{equation*}
Then,
\begin{equation*}
[L(z),L(w)]=\d_w L(w)\delta(z-w)+2L(w)\d_w\delta(z-w).
\end{equation*}

An element of a conformal superalgebra satisfying the above equation
is called a \emph{Virasoro element}; it is automatically
even.
\end{example}

In all other examples below we forgo the description of conformal
superalgebras in terms of formal distributions and simply provide the
(conformal) generators and relations.  A more detailed description can be
found in \cite{FK}.

\begin{example}\lbb{wser}  Recall that for the Grassmann algebra 
$\wedge(N)$ 
in the anti-commuting indeterminates $\xi_i$, $i=1,\dots,N$, its Lie 
superalgebra of derivations is
\begin{equation*}
W(N)=\left\{\sum_{i=1}^N P_i\d_i\st P_i\in\wedge(N), 
\d_i=\d/\d\xi_i\right\}.
\end{equation*}

The Lie conformal superalgebra $W_N$ is defined as
\begin{equation*}
W_N=\Cset[\d]\otimes \left( W(N) \oplus \wedge(N) \right).
\end{equation*}
The  $\la$-bracket $(a,b \in W(N); f,g \in \wedge(N))$ is as follows:
\begin{equation}\lbb{lbrwn}
[a_{\la}b]=[a,b],\;\;\; [a_{\la}f]=a(f)-\la (-1)^{p(a)p(f)}fa, \;\;\;
[f_{\la}g]=-\d(fg) -2\la fg.
\end{equation}

The Lie conformal superalgebra $W_N$ is simple for $N \geq 0$ and has rank 
$(N+1)2^N$.
We also remark that $W_0\simeq \vir$.

We shall need the following representation of $W_N$ on $\Cset[\d]\otimes
\wedge(N)$ (cf. Example~\ref{vironcur}):
\begin{equation}\lbb{lactwn}
a_{\la}g=a(g), \quad
f_{\la}g=-(\d + \la)fg, \quad  a\in W(N); f,g \in \wedge(N).
\end{equation}

As a consequence, by identifying $(\cur\g)\otimes\wedge(N)$ and
$\g\otimes\Cset[\d]\otimes\wedge(N)$, we also get a representation of
$W_N$ on $\cur\g\otimes\wedge(N)$ by conformal derivations.  Thus, we can 
construct the
semi-direct product $W_N\ltimes (\cur\g\otimes\wedge(N))$.

Recall that $W(N)$ and $\wedge(N)$ are graded, $W(N)=\oplus_{j\geq -1} 
W(N)^j$ and $\wedge(N)=\oplus_{j\geq 0} \wedge(N)^j$, with $\deg \xi_i=1$, 
$\deg 
\d_i=-1$.  A subalgebra $L$ of 
$W(N)$ is said to \emph{act transitively} on $\wedge(N)$ or 
$\g\otimes\wedge(N)$ if under the 
projection $W(N)\to W(N)^{-1}$, $L$ maps onto $W(N)^{-1}$.

Similarly, a subalgebra $L$ of $W_N$ \emph{acts transitively} on 
$\cur\g\otimes\wedge(N)$ if the projection of $L$ to $\CC[\d]\otimes 
W(N)^{-1}$ has rank $N$.
\end{example}

\begin{example}\lbb{sser}
For
an element $D=\sum_{i=1}^N P_i(\d,\xi)\d_i + f(\d,\xi) \in W_N$, we
define the corresponding notion of divergence:
\begin{equation*}
\Div D=\sum_{i=1}^N(-1)^{p(P_i)}\d_i P_i -\d f \in \Cset[\d]\otimes
\wedge(N).
\end{equation*}
The following identity holds in $\Cset[\d]\otimes
\wedge(N)$, where $D_1,D_2 \in W_N$ (cf. (\ref{lactwn})):
\begin{equation}\lbb{frmdiv}
\Div[{D_1}_{\la}D_2]=({D_1})_{\la}(\Div
D_2)-(-1)^{p(D_1)p(D_2)}({D_2})_{-\la-\d} (\Div D_1).
\end{equation}
Therefore we can define the following subalgebra of $W_N$:
\begin{equation*}
S_N=\{ D\in W_N \st \Div D=0 \}
\end{equation*}

The Lie conformal superalgebra $S_N$ is
simple for $N\geq 2$ and has rank $N2^N$.
\end{example}

\begin{example}\lbb{s2except2}
Let $D=\sum_{i=1}^N P_i(\d,\xi)\d_i + f(\d,\xi)$ be an
element of $W_N$.  Given $a\in\CC$, we define the deformed divergence to 
be
\begin{equation*}
\Div_a D= \Div D + a f.
\end{equation*}
It still satisfies formula
(\ref{frmdiv}), hence
\begin{equation*}
S_{N,a}= \{ D \in W_N \st \Div_a D = 0 \}
\end{equation*}
is a subalgebra of $W_N$, which is simple for $N \geq 2$ and has rank
$N2^N$.
\end{example}

\begin{example}\lbb{s2except1}  Another variation on the construction of the
conformal superalgebra $S_N$ is the following definition ($N$ even):
\begin{equation*}
\wti{S}_N=\{D \in W_N \st \Div((1+\xi_1 \ldots \xi_N)D)=0 \}
           (=(1-\xi_1 \ldots \xi_N)S_N).
\end{equation*}
We thus obtain the Lie conformal superalgebra $\wti{S}_N$. It is simple for $N
\geq 2$ and has rank $N2^N$.
\end{example}

\begin{example}\lbb{kser} We can also define the ``contact'' conformal 
superalgebra
\begin{equation*}
K_N=\Cset[\d] \otimes \wedge(N)
\end{equation*}
with the $\la$-bracket for $A=\xi_{i_1}\dots\xi_{i_r}$ and 
$B=\xi_{j_1}\dots\xi_{j_s}$ defined as
\begin{equation*}
[A_{\la}B]=\left(\left(\frac{r}{2}-1\right)\d (AB) + (-1)^r
\frac{1}{2}\sum_{i=1}^N
\d_iA\d_iB\right) +\la\left(\frac{r+s}{2}-2\right)AB.
\end{equation*}

$K_N$ embeds into $W_N$ but we will not be using this fact here.

Note also that $K_0\simeq\vir$ and $K_2\simeq W_1$.

The Lie conformal superalgebra $K_N$ is simple for all $N \in \Zset_+$,
$N \neq 4$ and is a free $\Cset[\d]$-module of rank $2^N$.
\end{example}

\begin{example}\lbb{k'4ex}  The Lie conformal superalgebra $K_4$ is not simple
but its derived subalgebra $K'_4$ is.
Furthermore,
$K_4=K'_4 \oplus \Cset \nu$, where $\nu=\xi_1\xi_2\xi_3\xi_4$ and,
since $K_4'$ is an ideal in $K_4$, the  map
$\ad \nu$ is an outer
conformal derivation of $K_4'$.
\end{example}

\begin{example}\lbb{ck6exa}
The Lie conformal superalgebra $CK_6$ is a simple rank 32
subalgebra of $K_6$, spanned over $\CC[\d]$ by the elements
\begin{equation*}
-1+\alpha\d^3\nu,\quad 
\xi_i-\alpha\d^2\xi_i^*,\quad\xi_i\xi_j-\alpha\d(\xi_i\xi_j)^*,
\quad \xi_i\xi_j\xi_k-\alpha(\xi_i\xi_j\xi_k)^*,
\end{equation*}
where $\alpha\in\CC$ is a fixed number such that $\alpha^2=-1$, 
$\nu=\xi_1\ldots\xi_6$, and 
$(\xi_{i_1}\xi_{i_2}\ldots)^*=\d_{i_1}\d_{i_2}\ldots\nu$.

The even part of $CK_6$ is $\vir \ltimes \cur
\so_6$ for the Virasoro element $-1+\alpha\d^3\nu$.
For the explicit form of the commutation relations of and further
details on $CK_6$, see \cite{CK2}.
\end{example}

The main result of \cite{FK} is the following Theorem.

\begin{theorem}\lbb{sthmfslcs}
Any finite simple Lie conformal superalgebra $R$ is isomorphic to one of
the Lie conformal superalgebras of the following list:
\begin{enumerate}
\item $W_N,$ $(N\geq 0)$;
\item $S_{N,a}$ $(N \geq 2, \;a \in \Cset)$;
\item $\wti{S}_N$ $(N \;even,\;\; N \geq 2)$;
\item $K_N$ $(N \geq 0, \;\; N\neq 4)$;
\item $K_4'$;
\item $CK_6$;
\item $\cur \ss$, where $\ss$ is a simple finite-dimensional Lie
superalgebra.
\end{enumerate}
\end{theorem}

We shall call the algebras (1)--(6) from the above list the Lie conformal 
superalgebras \emph{of Cartan type}.

\subsection{Conformal derivations of simple Lie conformal
superalgebras}\label{subsec.dersim}  We are going describe the
conformal derivations of simple Lie conformal superalgebras.
With this in mind, for every finite simple Lie conformal superalgebra
$R$ we will fix a \emph{distinguished} reductive Lie subalgebra $\rr$ of
$(R,\encirc{0})$.

For $R=\cur\ss$, we choose as $\rr$ the maximal reductive
subalgebra of $\ss$ (they are listed in \cite{K3}).

For $W_N$ we take as $\rr$ the 
copy of $\gl_N$ spanned by $\xi_i\d_j\in W(N)$.  We construct other 
distinguished subalgebras in a similar fashion.  For $S$-type 
superalgebras, we take the subalgebras $\sl_N$ of $\gl_N\subset W(N)$ and 
for 
the $K$-type, we take $\rr=\so_N$ spanned by $\xi_i\xi_j$.  For $K'_4$, 
$\rr=\so_4\oplus\CC\d\nu$ and for $CK_6$, 
the algebra $\so_6$ contained in $\cur\so_6$.  Summing up, we 
obtain the following list of $\rr$'s:
\begin{equation}\label{distlist}
\begin{split}
W_N,\; N \geq 0 \quad : \quad &\gl_N \\
S_{N,a},\; N \geq 2, \;  a\in \Cset \quad : \quad &\sl_N  \\
\wti{S}_N,\; N \geq 2, \; N\; \text{even} \quad : \quad &\sl_N \\
K_N,\;  N \geq 0,\; N \neq 4  \quad : \quad &\so_N \\
K'_4 \quad : \quad &\cso_4 \\
CK_6 \quad : \quad &\so_6
\end{split}
\end{equation}

\begin{proposition}\lbb{cderfslcs}
In the following we list all cases in which a finite simple Lie
conformal superalgebra has outer conformal derivations.
\begin{enumerate}
\item $\cder K_4'= \cinder K_4' \oplus \Cset[\d]\nu$,
where $\nu=\xi_1\dots\xi_4\in\wedge(4)$;
\item $\cder (\cur  \ss) =\vir \ltimes \cur (\der \ss)$,
where $\ss$ is a simple finite-dimensional Lie superalgebra.
\end{enumerate}
\end{proposition}

\begin{proof} Let $R$ be a
finite simple Lie conformal superalgebra.  Then $R=\CC[\d]\otimes U$ and 
the distinguished reductive Lie algebra $\rr$ acts via the $0$-th product
completely reducibly on $U$.  This action commutes with $\d$, hence $\rr$ 
acts completely reducibly on $R$.

Next, $\cder R\subset\gc R\simeq R\otimes R^*$ as an $R$-module, where 
$R^*$ is the conformal dual of the $R$-module $R$ \cite[Proposition 
6.4(a)]{BKL}.  Therefore $\gc R$ and, hence, $\cder 
R$ are completely reducible as $\rr$-modules.

Thus, $\cder R= \cinder R\oplus V$ as $\rr$-modules. We will refer to the
elements in $V$ as \emph{outer} conformal derivations of $R$.

Since $[\rr\encirc{0} V]\subset R=\cinder R$, we see that $\rr$ kills $V$,
thus outer conformal derivations are $\rr$-module homomorphisms.
Using this, outer conformal derivations of $R$ can be easily computed.

In the following we provide a detailed computation in the most
involved case, that of $W_N$. We restrict our attention to $N \geq
2$, the remaining cases being straightforward. 

Recall that $W_N=\Cset[\d]\otimes (W(N) \oplus \wedge(N))$ and that $W_N$ 
is completely reducible as a $\gl_N$-module.  We have the following 
description of $W(N)^k\oplus\wedge(N)^k$ as $\sl_N$($\subset\rr$)-modules.  
Here $R(\lambda)$ denotes the $\sl_N$-module with the the highest weight 
$\lambda$ and  $\pi_i$ are fundamental weights:
\begin{equation*}
W(N)^k\oplus\wedge(N)^k=
\begin{cases}
R(\pi_{N-1}), \quad\text{if } k=-1,\\
R(\pi_{k+1} +\pi_{N-1}) \oplus R(\pi_k) \oplus R(\pi_k),\quad\text{if } 
0 \leq k \leq N-2,\\
R(\pi_{N-1})\oplus R(\pi_{N-1}),\quad\text{if } k= N-1,\\
R(0),\quad\text{if }k=N. 
\end{cases}
\end{equation*}
Let $\phi \in V$. Then $\phi$ is a $\gl_N$-module homomorphism, hence it
commutes with the Euler operator $E=\sum\xi_i\d_i$ and therefore leaves
each graded component invariant. The grading being consistent with parity, 
$\phi$ can only be an even conformal derivation.
Let $x_k=\xi_1\ldots\xi_k E$ be the highest weight vector of $R(\pi_k)
\subseteq W(N)^k$, $0 \leq k \leq N-1$.
Let $z_k=\xi_1 \ldots \xi_{k+1} \d_N$ be the highest weight vector of
$R(\pi_{k+1} + \pi_{N-1})\subseteq W(N)^k$, $0 \leq k \leq N-2$. Let
$\d_N$
be the highest weight vector of $R(\pi_{N-1}) \subseteq W(N)^{-1}$.
Let $y_k=\xi_1 \ldots \xi_k$  be the highest weight vector of $R(\pi_k)
\subseteq \wedge(N)^k$, $0 \leq k \leq N$. The action of $\phi$ is
given by
\begin{align*}
&\phi_{\la}y_N=A(\d,\la)y_N,  \\
&\phi_{\la}x_k=P_k(\d,\la)x_k + Q_k(\d,\la)y_k,  \quad 0 \leq k\leq N-1,\\
&\phi_{\la}y_k=R_k(\d,\la)x_k + S_k(\d,\la)y_k, \quad 0 \leq k\leq N-1,\\
&\phi_{\la}z_k=Z_k(\d,\la)z_k,  \quad 0 \leq  k \leq N-2,\\
&\phi_{\la}\d_N =\Omega(\d,\la)\d_N.
\end{align*}
Let $g\in \sl_N$. We compute $\phi_{\la}[g_{\mu}x_k]$: 
\begin{align*}
\phi_{\la}[g_{\mu}x_k]=&(\phi_\la g)_{\la+\mu}x_k+g_\mu(P_k(\d,\la)x_k + 
Q_k(\d,\la)y_k)\\
=&(\phi_\la g)_{\la+\mu}x_k+P_k(\d+\mu,\la)(g_\mu 
x_k)+Q_k(\d+\mu,\la)(g_\mu 
y_k).
\end{align*}
Since for any $g\in\sl_N$, there exists $g'$ such that $g=[g',z_0]$, we 
have $\phi_\la g=\phi_\la(g'\encirc{0}z_0)=Z_0(\d,\la)g$ as $\phi_\la$ and 
$g'\encirc{0}$ commute.  Then using (\ref{lbrwn}) we get
\begin{equation*}
\phi_\la[g,x_k]=
Z_0(-\la-\mu,\la)[g,x_k]+P_k(\d+\mu,\la)[g,x_k]+Q_k(\d+\mu,\la)(gy_k-\mu y_kg).
\end{equation*}
By commuting $\phi_\la$ and $g\encirc{0}$, we get that the left-hand side 
of the above equality equals $P_k(\d,\la)[g,x_k]+Q_k(\d,\la)gy_k$. 

There exists $g\in\sl_N$ such that $y_kg\neq [g,x_k]$ (e.g. 
$g=\xi_1\d_1$).  Hence, 
$Q_k(\d+\mu,\la)=0$ implying that $Q_k=0$.  Also, by comparing the 
left-hand and the right-hand sides, we see that $P_k(\d,\la)$ is at most 
linear in $\d$, otherwise it is impossible to cancel products 
of $\mu$ and $\d$ that may come from $P_k(\d+\mu,\la)$. 

The rest of the deduction is similar and we describe it only in brief.
We compute $\phi_{\la}[{x_0}_{\mu} x_k]$,
$\phi_{\la}[{y_0}_{\mu} y_k]$ and $\phi_{\la}[g_{\mu}y_k]$ to show that
$R_k(\d,\la)$ is constant in $\d$ whereas $Z_k(\d, \la)$ and
$S_k(\d,\la)$ are at most linear in $\d$. Third, we compute
$\phi_{\la} [{x_0}_{\mu}y_k]$,
$\phi_{\la} [{x_0}_{\mu}\d_N]$,
$\phi_{\la} [{\d_N}_{\mu}y_0]$,
$\phi_{\la} [{x_k}_{\mu}y_0]$,
$\phi_{\la} [{ x_0 }_{\mu}y_N]$ and
$\phi_{\la} [{\d_N}_{\mu}\xi_N]$ to show that all the remaining
polynomials $A, P_k, R_k, S_k, Z_k, \Omega_k$ are actually zero.
\end{proof}

\begin{remark} Using the same arguments as in the proof of 
Proposition~\ref{cderfslcs}, one can calculate the ordinary derivations of 
finite simple Lie conformal superalgebras.  Below we list all cases in 
which 
outer ordinary derivations occur ($E$ denotes the Euler operator and 
$\nu$, the highest monomial in $\wedge(N)$):
\begin{enumerate}
\item
$\der S_{N,a}=\inder S_{N,a} \oplus \Cset E_{(0)} \;\;
a\neq 0, N\geq 2$;
\item
$\der S_{N,0} =\inder S_{N,0} \oplus \Cset E_{(0)} \oplus
\Cset \nu_{(0)}, \;\;N >2$;
\item
$\der S_{2,0} =\inder S_{2,0} \oplus \sl_2$, where 
$E\encirc{0},\nu\encirc{0}\in\sl_2$;
\item
$\der \wti{S}_N =\inder \wti{S}_N \oplus
\Cset \nu_{(0)} \;\;N \geq 2$, $N$ even;
\item
$\der K_4'=\inder K_4'  \oplus \Cset \nu_{(0)}$;
\item
$\der \cur \ss =\Cset \d \oplus \der \ss$, where $\ss$ is a simple
finite-dimensional Lie superalgebra.
\end{enumerate}

We will not require this result for the rest of the paper and so leave the 
details to the reader.

It follows that the ordinary derivations of simple finite Lie conformal 
superalgebras form the following Lie superalgebras:
\begin{enumerate}
\item $\der W_N\simeq W(N)\ltimes \wedge(N)$;
\item $\der S_{N,a}\simeq (S(N)\oplus\CC E)\ltimes\wedge(N)'$, where 
$\wedge(N)'$ is the subalgebra of $\wedge(N)$ spanned by monomials 
of degree strictly less than $N$, $a\neq 0$;
\item      $\der S_{N,0}\simeq (S(N)\oplus\CC E)\ltimes\wedge(N)$, 
$N>2$;
\item      $\der S_{2,0}\simeq \so_4\ltimes \mathcal{H}_4$, where 
$\mathcal{H}_4$ is the Heisenberg algebra generated by 
$\xi_1,\xi_2,\d_1,\d_2$, and $\so_4$ is spanned by two $\sl_2$-triples: 
$(\xi_2\d_1, \xi_1\d_1-\xi_2\d_2,\xi_1\d_2)$, and $(\nu, E, F)$, where 
$F(\xi_1)=-\d_2$, $F(\xi_2)=\d_1$.
\item $\der\wti{S}_N\simeq ((1-\nu)S(N))\ltimes\wedge(N)$, $N$ even;
\item $\der K_N$, $N\neq 4$, is the natural central extension of $H(N)$ 
by a one-dimensional center;
\item $\der K'_4\simeq H(4)$;
\item $\der CK_6$ is the central extension of the ``strange" Lie 
superalgebra $P(4)$ by a one-dimensional center;
\item $\der \cur\ss$ is the direct sum of $\der\ss$ with the 
one-dimensional Lie algebra.
\end{enumerate}
(Cf. the calculations of derivations of corresponding formal distribution 
Lie superalgebras in \cite{K5} and \cite{FK}.)
\end{remark}


\section{Conformal derivations of differentiably simple Lie conformal
superalgebras}\label{sec.diffsimder}

Here we describe conformal derivations of Lie conformal superalgebras of
the form $S\otimes\wedge(n)$, where $S$ is a finite simple Lie conformal
superalgebra.  By Theorem \ref{diffsimple} this will take care of the
conformal derivations of all differentiably simple Lie conformal
superalgebras.

\subsection{Conformal centroid} The {\em conformal centroid} of a Lie
conformal superalgebra $R$ is the subalgebra $CC(R)$ of
the associative conformal algebra $\cend R$ defined by
\begin{equation*}
CC(R)=\{\varphi\in\cend R\st \varphi_\lambda [x_\mu y]=[(\varphi_\lambda 
x)_{\la+\mu} y]\}.
\end{equation*}

\begin{remark}\label{ccr1} A direct calculation similar to that in the 
proof of Lemma~\ref{commcent} shows that for $\varphi\in CC(R)$,
$\varphi_\lambda [x_\mu y]=(-1)^{p(x)p(\varphi)} [x_\mu(\varphi_\lambda
y)]$.
\end{remark}

To describe conformal centroids of simple superalgebras, we need two
technical lemmas:

\begin{lemma}\lbb{lmcoce}
Let $R$ be a finite Lie conformal superalgebra.
Suppose $R$ contains a Virasoro element $L$
and a $\Cset[\d]$-basis $\{ a_j \}_{j \in J}$ of $R$ such that
\begin{equation}\lbb{eqvir}
[L_{\la} a_j]=(\d + \Delta_j \la)a_j + b_j \quad \text{for any} \quad
j \in J,
\end{equation}
where $\Delta_j
\in \Cset$, $b_j\in R$ and, whenever $\Delta_j =0$, $b_j=0$.

Then $CC(R)=0$.
\end{lemma}

\begin{proof}
Let $\varphi \in CC(R)$. Suppose $\deg_\lambda \varphi_{\la} L=k$.
Then $\deg_\lambda [L_\mu (\varphi_\lambda L)]=k$ and $\deg_\lambda
\varphi_\lambda [L_\mu L]=\deg_\lambda (\d+\lambda+2\mu) \varphi_\lambda
L=k+1$ unless $\varphi_{\la} L=0$.  Hence the equation
$\varphi_{\la}[L_{\mu}L]=[L_{\mu}(\varphi_{\la}L)]$, given in 
Remark~\ref{ccr1} implies that $\varphi_{\la} L=0$.

Now, for any $j \in J$, we have
\begin{equation*}
\varphi_{\la}[L_{\mu}a_j]=\varphi_{\la}((\d +\Delta_j \mu)a_j + d(a_j)).
\end{equation*}
On the other hand,
\begin{equation*}
\varphi_{\la}[L_{\mu}a_j]=[(\varphi_{\la}L)_{\la + \mu} a_j]=0.
\end{equation*}
Therefore,
\begin{equation*}
(\d + \la + \Delta_j \mu)\varphi_{\la}a_j + \varphi_{\la}b_j =0.
\end{equation*}
Now, if $\Delta_j \neq 0$ the fact that there is no power of $\mu$ in the
second summand implies that $\varphi_{\la}a_j=0$. If $\Delta_j = 0$,
then,
by assumption, $b_j=0$, hence  $\varphi_{\la}a_j=0$.
\end{proof}

Now we present a method of constructing a basis satisfying requirements of 
Lemma~\ref{lmcoce} which we will apply in the case of 
simple Lie conformal superalgebras of Cartan type.

\begin{lemma}\lbb{exvire}
Let $R$ be a finite simple Lie conformal superalgebra, $R=\CC[\d]\otimes 
V$.  Let $\rr$ be the distinguished reductive Lie subalgebra of $R$ and 
let $L$ be a Virasoro element of $R$ such that 
\begin{equation}\label{eqvir2}
[L_\la g]=(\d+\la)g\quad \text{for all }g\in\rr.
\end{equation} 
Suppose $V=\oplus_{i \in I} V_i$, where $V_i$ is an 
irreducible $\rr$-module whose highest weight vector is $v_i$. Suppose 
that the $v_i$ satisfy (\ref{eqvir}).
Then
there exists a $\Cset[\d]$-basis $\{ a_j \}_{j \in J}$ of $R$ satisfying 
conditions of Lemma~\ref{lmcoce}.
\end{lemma}

\begin{proof}
For any $i\in I$, $V_i$ is spanned over $\Cset$ by  vectors of the form
$v_i^n=g^n_{(0)}\ldots g^1_{(0)}v_i$, $n \in \Zset_+$ and $g^k \in \rr$
for any $k=1, \ldots,n$.

We are going to prove by induction on $n$ that
\begin{equation*}
[L_{\la}v_i^n]=(\d + \Delta_i\la)v_i^n + u_i^n,
\end{equation*}
where $u_i^n\in R$, and whenever $\Delta_i=0$, $u_i^n=0$ for all $n$.
Then it would suffice to take $v_i^n$'s as the necessary basis.  We have
\begin{align*}
[L_{\la}(g^{n+1}_{(0)}v_i^n)]&=[[L_{\la}g^{n+1}]_{\la}v_i^n] + 
[g^{n+1}_{(0)}[L_{\la}v_i^n]]
\\
&=[((\d +\la)g^{n+1})_{\la}v_i^n] + [g^{n+1}_{(0)}((\d +\Delta_i \la)v_i^n 
+ u_i^n]
\\
&= (\d +\Delta_i \la)[g^{n+1}_{(0)}v_i^n] +[g^{n+1}_{(0)}u_i^n].
\end{align*}
Then we take $u_i^{n+1}=[g^{n+1}_{(0)}u_i^n]$.  This completes the 
proof.
\end{proof}

\begin{proposition}\lbb{cocefslcs}
Let $R$ be a finite simple Lie conformal superalgebra.
\begin{enumerate}
\item If $R$ is of Cartan type, then $CC(R)=0$;
\item If $R=\cur\ss$, then $CC(R)=\{\varphi\in\cend R\st \varphi_\lambda 
s=p(\d,\lambda)s, 
p(\d,\la)\in\CC[\d,\la], \text{for all }s\in\ss\}$.
\end{enumerate}
\end{proposition}

\begin{proof}
If $R$ is of Cartan type, then it has an element $L$ and a
$\Cset[\d]$-basis satisfying (\ref{eqvir}) and (\ref{eqvir2}) (such an $L$ 
is called physical in Section~\ref{sec.physpairs}). In fact, if $R$ 
is not
$S_{N,a}$ or $\wti{S}_N$, we can choose $L=-1$ (and then all $b_j=0$).
In the case of $S_{N,a}$ or $\wti{S}_N$, we choose 
$L=-1+\frac{1}{N}(\d-a)E$ and $-1+\xi_1\dots\xi_N+\frac{1}{N}\d E$, 
respectively. The formulas 
that appear in the proof of \cite[Proposition 4.16]{FK} show that
the highest weight vectors satisfy (\ref{eqvir}) with 
$\Delta_i$ taking the values of $2/N, 3/N,
\ldots, (N-1)/N, 1, (N+1)/N, \ldots ,(2N -1)/N, 2$. 
Consequently, by 
Lemma \ref{exvire} we can find a $\Cset[\d]$-basis whose conformal weights 
are nonzero. We complete the proof of (1) by applying Lemma~\ref{lmcoce}. 

Let $\varphi \in CC(\cur \ss)$. Then for any $s,s' \in \ss$ we have
\begin{equation*}
\varphi_{\la}[s_{(0)}s']=(-1)^{p(\varphi)p(s)}[s_{(0)}\varphi_{\la}s'].
\end{equation*}
Suppose $\varphi_{\la}=\sum_{n \in \Zset_+}  \frac{\la^n}{n!}\varphi_{(n)}$.
Then for any
$n \geq 0$, $\varphi_{(n)}$ is a $\ss$-module homomorphism.
Since $\ss$ is simple, if $\varphi_{(n)}\neq 0$, then
$\ker \varphi_{(n)}=0$ and $\im \varphi_{(n)}\simeq \ss$.  Then 
$\phi\encirc{n}s=\sum_i p_{n_i}(\d) s_i(s)$.  The map $s\mapsto s_i(s)$ is 
an 
$\ss$-automorphism.  Remark~\ref{schurlie} implies that 
$s_i(s)=c_is$ for some $c_i\in\CC$.  Hence, $\varphi_{(n)}(s)= 
p_n(\d) \otimes s$ for any $s \in \ss$ and
$\varphi_{\la}= \sum_{n \in \Zset_+} \frac{\la^n}{n !}p_n(\d)
\otimes 1_{\ss}.$
\end{proof}

\begin{remark} In fact, we have shown that $CC(\cur\ss)\simeq\cend_1$.  We 
note that this conformal algebra is neither finite, nor commutative.
\end{remark}

\subsection{Conformal derivations of a tensor product}

\begin{proposition}\lbb{cdertp}
Let $R$ be a finite Lie conformal superalgebra and let $B$
be a unital commutative associative finite-dimensional superalgebra. Then
\begin{enumerate}
\item
$\cder R \otimes B \subseteq \cder(R\otimes  B)$.
\item
$CC(R)\otimes \der B \subseteq \cder(R\otimes B).$
\item If $R$ is simple and of Cartan type, then
\begin{equation*}
\cder(R \otimes B ) =  \cder R\otimes  B.
\end{equation*}
\item If $R=\cur \ss$, where $\ss$ is a simple finite-dimensional Lie
superalgebra, then
\begin{equation*}
\cder( \cur  \ss \otimes B) = \cder(\cur \ss) \otimes B +  1_{\ss}
\otimes
\cder \cur  B.
\end{equation*}
In particular,
\begin{equation*}
\cder( \cur  \ss \otimes \wedge(N)) = \cur(\der \ss) \otimes
\wedge(N) + 1_\ss\otimes W_N.
\end{equation*}
\end{enumerate}
\end{proposition}

\begin{proof}
(1) Let $\phi\in \cder R$ and $b \in B$.
We set
\begin{equation*}
(\phi\otimes b)_{\la}(r' \otimes b')=(-1)^{p(b)p(r')}(\phi_{\la}r')
\otimes bb'.
\end{equation*}
It is easy to verify that $\phi \otimes b$ is a conformal derivation.

(2) Similarly, let $\varphi \in CC(R)$ and $d \in \der B$. We set
\begin{equation*}
(\varphi \otimes d)_{\la}(r \otimes b)=(-1)^{p(d)p(r)}(\varphi_{\la}r)
\otimes
d(b).
\end{equation*}
It is easy to verify that $\varphi \otimes d$ is a conformal
derivation.

(3) Let $\phi \in \cder(R \otimes B)$. For any $r \in R$ we have
\begin{equation*}
\phi_{\la}(r \otimes 1)=\sum_{i \in I} {\phi_i}_{\la}(r) \otimes b_i
\end{equation*}
where $\{b_i\}_{i \in I}$ is a linear basis of $B$. A short 
direct computation
shows
that $\phi_i \in \cder R$ and $\wti{\phi}:=\phi-\sum_{i \in I} \phi_i
\otimes b_i \in
\cder(R \otimes B)$ is such that $\wti{\phi}_{\la}(r \otimes 1)=0$.
Let us fix $b \in B$.
Suppose $\wti{\phi}_{\la}(r \otimes b)=\sum_{i \in I} {\varphi_i}_{\la}r
\otimes b_i$.
For any $r,r' \in R$ we have
$[r_{\mu}r']\otimes b=[(r\otimes 1)_{\mu}(r' \otimes b)]$.
If we apply $\wti{\phi}$ to both sides of this equation we see that
$\varphi_i
\in CC(R)$ which is zero by Proposition \ref{cocefslcs}.

(4) Using the same argument as in (3), one can see that $\varphi_i \in
CC(\cur \ss)$, so that by Proposition \ref{cocefslcs}
\begin{equation*}
\wti{\phi}_{\la}(r \otimes b)=\sum_{i \in I} P_i(\d, \la) r \otimes b_i.
\end{equation*}
Now, we identify $\cur \ss \otimes B$ with $\ss \otimes \cur  B$ so that
\begin{equation*}
\wti{\phi}_{\la}(r \otimes b)=r \otimes \sum_{i \in I} P_i(\d, \la) 
\otimes b_i.
\end{equation*}
The map associating $b$ to $\sum_{i \in I} P_i(\d, \la) \otimes b_i$
is easily shown to be a conformal derivation of $\cur B$.

In the case $B=\wedge(N)$, it remains to show that $W_N= \cder(\cur \wedge 
(N))$. We recall at first that   $W_N$ acts on $\cur \wedge(N)$ by 
conformal derivations (see (\ref{lactwn})).
The rest is done  in two steps. First, notice that if
$\phi \in \cder(\cur \wedge (N))$, then $\phi_\la 1 =(\d + \la)
\sum_{i \in I} p_i(\la)f_i$,
where $\{f_i\}_{i \in I}$ is a linear basis of $\wedge(N)$. Second,
if we set $\wti{\phi} =\phi -(\sum_{ i \in I} p_i(-\d)f_i) \circ \d$,
it is immediate to see that
for any $j=1, \ldots , N$, $\wti{\phi}_{\la}(\xi_j)= P_j(\la, \xi)$,
so that $\wti{\phi}=\sum_{j=1}^N
P_j(-\d, \xi) \d_j$. Therefore $\phi= \wti{\phi} + (\sum_{i \in I}
p_i(-\d) f_i) \circ \d \in W_N$.
\end{proof}

We thus obtain

\begin{theorem}\label{diffsimple1} The following is a complete
list of finite non-abelian differentiably simple Lie conformal 
superalgebras:
\begin{enumerate}
\item $(\cur\ss)\otimes\wedge(n)$, where $\ss$ is a simple
finite-dimensional Lie superalgebra;
\item a finite simple Lie conformal superalgebra of Cartan type.
\end{enumerate}
\end{theorem}

\begin{proof} If $R$ is differentiably simple, then $R\simeq
S\otimes\wedge(n)$ by Theorem~\ref{diffsimple}.  But if $S$ is
of Cartan type, an ideal $S\otimes I$ of $R$, where $I$ is an ideal
of $\wedge(n)$, is differentiably stable by
Proposition~\ref{cdertp}(3).  Part (4) of the same proposition
shows that the conformal superalgebra $(\cur\ss)\otimes\wedge(n)$,
$\ss$ as above, is differentiably simple.
\end{proof}


\section{Semisimple Lie conformal
superalgebras}\label{sec.semisim}

Here we provide a detailed description of finite semisimple Lie conformal
superalgebras.  In particular, we prove the following

\begin{theorem}\label{mainth} Let $R$ be a  finite semisimple Lie conformal
superalgebra.  Then $R$ splits into a direct sum of conformal algebras of
the following types:
\begin{enumerate}
\item a finite simple Lie conformal superalgebra;
\item a Lie conformal superalgebra $L$ such that
\begin{multline}\label{eq.mainth}
\left(\bigoplus_{i=1}^k (\cur \ss_i)\otimes\wedge(n_i)\right)\oplus
{K'_4}^{\oplus^r}
\subset L\\ \subset
\left(\bigoplus_{i=1}^k \left( 
\cur(\der\ss_i)\otimes\wedge(n_i)+1_\ss\otimes 
W_{n_i}\right)\right)\oplus
{K_4}^{\oplus^r},
\end{multline}
where $n_i,r\in\ZZ_+$ and $\ss_i$ are simple finite-dimensional Lie 
superalgebras,
and such that for each $i$ the projection of $L$ onto $W_{n_i}$ acts 
transitively on $(\cur\ss_i)\otimes\wedge(n_i)$.
\end{enumerate}
\end{theorem}

\subsection{Proof of Theorem~\ref{mainth}}  Define the {\em
socle} of a Lie conformal algebras as the sum of all its minimal ideals.

\begin{lemma}\lbb{smiisdi}
Let $R$ be a finite semisimple Lie conformal superalgebra. Then
its socle is a direct sum of all minimal ideals of $R$ and
there exist finitely many such ideals in $R$.
\end{lemma}

\begin{proof}
Let $\{M_i \}_{i \in I}$ be the family of all minimal ideal of $R$.
By Lemma \ref{milcs} this family is non-empty.

Notice that for all $i,j \in I$, $[{M_i}_\lambda M_j]=0$ for $i\neq j$ by
minimality.  Suppose $M_i \cap \sum_{j \neq i} M_j \neq 0$ for some $i$.
Thus by minimality of $M_i$,  $M_i \subseteq \sum_{j \neq i}
M_j$, and $[{M_i}_\lambda M_i]\subset \sum_{j \neq i} [{M_i}_\lambda M_j]
=0$.  Hence $M_i$ is abelian which contradicts semisimplicity of $R$.

Thus the sum of minimal ideals is direct. Furthermore, the rank of a
direct sum is the sum  of  the ranks, hence the fact that $R$ is finite
implies that $R$ contains only finitely many minimal ideals.
\end{proof}

We fix notations from the proof and denote the socle of $R$ by $M$ and the
minimal ideals of $R$ by $M_i$, $i=0,1,\dots,l$.

The centralizer of $M$ in $R$, $C_R(M)=\{x\in R \st [x_{\la} M]=0\}$, is
zero.  This follows from the fact that $C_R(M)$ is an ideal of $R$, but
contains no minimal ideals of $R$, because any such ideal would be abelian
contradicting the semisimplicity of $R$.  But this contradicts Lemma
\ref{milcs}.

Therefore, the homomorphism of Lie conformal superalgebras
$R \to \cder M$ sending $x$ to $\ad_M x_{\la}$ is injective and we have
\begin{equation*}
\cinder M \subseteq R \subseteq \cder M.
\end{equation*}

The minimality of $M_i$ implies that $\langle [{M_i}_{\la} M_i] \rangle =
M_i$ for any $i=0,1,\ldots, l$.  Thus if $\phi\in\cder M$ maps $M_i$ to
$M_j$, $\phi(M_i)$ must be zero.  It follows that $\cder M=\oplus_{i=0}^l
\cder M_i$.

Therefore, we have
\begin{equation}\label{eq.inbetween}
\bigoplus_{i=0}^l \cinder M_i \subseteq R \subseteq \bigoplus_{i=0}^l
\cder M_i.
\end{equation}

Put $R_i=R\cap \cder M_i$.

\begin{lemma}\label{mainsemisim} Let $R$ be a conformal superalgebra 
satisfying (\ref{eq.inbetween}).  Then $R$ is semisimple if and only if 
$M_i$ is $R_i$-simple for all $i$.
\end{lemma}

\begin{proof} If $M_i$ is not $R_i$-simple for some $i$, then it contains
a non-trivial ideal $J$. By differential simplicity of $M_i$, $J$ is
nilpotent (see Proposition~\ref{maxid}).  Since $J$ is also an ideal of
$R_i$ and hence of $R$, $R$ is not semisimple.

Conversely, let $R$ contain an abelian ideal $I$.  Since $M_i$ is 
non-abelian and $R$-simple for any $i$, $[I_\la M_i]$, as an 
$R$-submodule of $M_i$ must be zero. Thus $I$ kills 
$M$, a contradiction. 
\end{proof}

By Lemma~\ref{fctsmi}(2), $M_i$'s are differentiably simple and, since 
they  are non-abelian, we obtain a complete description of the possible 
form of $M_i$'s from Theorem~\ref{diffsimple1}.

Now we can use Propositions \ref{cderfslcs} and \ref{cdertp} to describe
possible $R_i$'s in greater detail.

Let $M_i=S_i\otimes\wedge(n_i)$, where $S_i$ is a simple Lie conformal
superalgebra.

By Proposition~\ref{diffsimple1}, if $S_i$ is of Cartan type, $n_i=0$.

If $S_i$ is of Cartan type and not isomorphic to $K'_4$, then $R_i\simeq
S_i$.

If $S_i\simeq K'_4$, $R_i\subset K_4$.

It remains to treat the case of a current Lie conformal algebra $S_i$.
Let $R_i$ be a semisimple Lie conformal superalgebra such that $\cur\ss 
\otimes\wedge(n)\subset R_i\subset W_n\ltimes 
(\cur\der\ss\otimes\wedge(n))$, 
where $\ss$ is a finite-dimensional simple Lie superalgebra.
In order to prove that the projection of $R_i$ onto $W_{n_i}$ acts 
transitively if and only if  $S_i\otimes\wedge(n_i)$ is $R_i$-simple, we
need an auxiliary lemma

Recall that a commutative associate superalgebra with the action of a Lie 
superalgebra $\fa$ by derivations is called $\fa$-differentiably simple if 
it contains no non-trivial invariant ideals.

\begin{lemma}\label{wnsubalg1}  $\wedge(n)$ is $\fa$-simple for a 
subalgebra  $\fa$ of $W(n)$ if and only if $\fa$ acts transitively on 
$\wedge(n)$.
\end{lemma}

\begin{proof} Transitivity of the action of $\fa$ is equivalent to saying 
that $\fa$ contains elements $a_j=\d_j+\sum_i p_i\d_i$, where 
$p_i\in\wedge(n)$, $p_i(0)=0$, for all $j$.

Then, assume that $\fa$ acts transitively and let $I$ be a non-zero 
$\fa$-stable ideal of $\wedge(n)$.  Then $I$ contains the monomial 
$\xi_1\dots\xi_n$.  Let $f=\xi_{k_1}\dots\xi_{k_n}$.  Since 
$a_j(\xi_jf)=f+$ monomials of higher degree, by induction $I$ 
contains all monomials of $\wedge(n)$, i.e. $I=\wedge(n)$.

Now let $\wedge(N)$ be $\fa$-differentiably simple.  We can assume that 
$\fa$ is closed
with respect to multiplication by $\wedge(N)$.  Indeed, $[a,fb]
=a(f)b+(-1)^{p(a)+p(f)}[a,b]$ for $a,b\in W(N)$, $a$ homogeneous, so
$\wedge(N)\fa$ remains closed with respect to the Lie bracket.  Also,
$\fa$ acts transitively if and only if $\wedge(N)\fa$ does.  Similarly, 
$\wedge(N)$ is $\fa$-differentiably simple if and only of it is 
$\wedge(N)$-differentiably simple.

If the projection of $\fa$ to $W(N)^{-1}$ is zero, then by the grading
argument, $\fa$ leaves invariant the maximal ideal of $\wedge(N)$.  Thus,
possibly after a linear transformation of $\xi_i$'s, $\fa$ contains an
element of the form $a_N=\d_N+\sum p_i\d_i$, $p_i(0)=0$.

In particular, $a_N=\d_N+\xi_1 b+\sum q_i\d_i$, where $b\in W(N)$,
$q_i(0)=0$ and $\d_1q_i=0$.  Also, we can assume that $q_N=0$. We
introduce the following change of coordinates: $\xi_i'=\xi_i-\xi_Nq_i$ for
$1\leq i\leq N-1$, $\xi_N'=\xi_N$.  Notice that
$\CC\<\xi_1',\xi_2\dots,\xi_N\>=\wedge(N)$.  Indeed, let $q_1=\xi_1f+g$,
$\d_1g=0$.  Then $\xi_1=\xi_1(1+\xi_Nf)+\xi_Ng$, and since $1-f$ is
invertible, this change of coordinates is valid.  By induction, we obtain
that $\CC\<\xi_1',\dots,\xi_N'\>$.  The chain rule $\d_i=\sum
(\d_i\xi_j')\d_j'$ implies that in the new coordinates $a_N=\d_N'+\xi_N c$
for some $c\in W(N)$.

Thus, we can assume that $a_N=\d_N+\xi_N c$. Then $[a_N,a_N]=2c$, hence 
$\d_N\in\fa$ and we 
can split $\fa$ as $\wedge(N)\d_N\oplus \fa_1$. If $\fa_1$ 
contains an element $a=\xi_N b$, then $[\d_N,a]=b\in\fa_1$, i.e. we can 
both multiply by and cancel $\xi_N$ in elements of $\fa_1$.  Thus $\fa_1$ 
splits as $\fa_2\oplus\xi_N\fa_2$, where $\fa_2\subset W(N-1)$.  If $I$ is 
a proper $\fa_2$-stable ideal of $\wedge(N-1)$, then clearly $I+\xi_1 I$ 
is a proper $\fa$-stable ideal.  Hence, $\wedge(N-1)$ is 
$\fa_2$-differentiably simple.  Moreover, $\fa$ acts transitively if and 
only if $\fa_2$ does.

Induction on $N$ completes the proof.
\end{proof}

\begin{remark} Along the way we proved an analog of the
Frobenius theorem for $\wedge(N)$.  Namely, if $I$ is an ideal of 
$\wedge(N)$ and $\fa$ is a subalgebra of $W(N)$ which leaves $I$ 
invariant, then after a change of coordinates, $\fa=\sum_{i=1}^k 
\wedge(N)\d_i+\sum_{j=k+1}^N \wedge_k(N)\d_j$, where $\wedge_k(N)$ is the 
subalgebra (without $1$) of $\wedge(N)$ generated by 
$\xi_{k+1},\dots,\xi_N$. 
\end{remark}

The explicit description of the action of $W_n$ on 
$\cur\ss\otimes\wedge(n)$ in the proof of Proposition~\ref{cdertp}(4) 
shows that the action of $R_i$ is transitive if and only if the action of 
$R_i|_{\d=1}$ is transitive on $\wedge(n)$. 
This together with Lemmas~\ref{mainsemisim} and~\ref{wnsubalg1} completes 
the proof of Theorem~\ref{mainth}.

\subsection{Conformal derivations} The following is a direct analog of
\cite[Proposition 7.4]{C}:

\begin{proposition} Let $R$ be a finite semisimple Lie conformal 
superalgebra and let $M$ be its socle.  Then $\cder R=\{\phi\in\cder M\st 
\<\phi_\lambda R\>\subset R\}$, the normalizer of $R$ in $\cder M$.
\end{proposition}

\begin{proof} If $I$ is a minimal $\cder R$-stable ideal,
it is either minimal or properly contains an ideal $J$ of $R$.  Since $I$
is differentiably simple, $J$ must be nilpotent, a contradiction.  Hence
$I\subset M$.  Conversely, let $M'$ be the sum of all minimal $\cder
R$-stable ideal and $I$ a minimal ideal such that $I\not\subset M'$.  Then
$I\subset\ann M'$ and, as $\ann M'$ is an ideal of $\cder
R$ and contains no $\cder R$-central elements, it must contain a minimal
$\cder R$-ideal which must be abelian.  The contradiction shows that
$M=M'$ and $M$ is $\cder R$-stable.

Then $\cder R$ embeds into $\cder M$ and the rest follows
easily.
\end{proof}


\section{Physical Virasoro pairs}\label{sec.physpairs}

In this section we compute all physical Virasoro pairs in a finite
simple Lie conformal superalgebra. As a consequence, we obtain the
classification of physical Lie conformal superalgebras.

\subsection{Definitions}

Let $R=\CC[\d]\otimes V$ be a finite simple Lie conformal superalgebra.
Recall that $V$ contains a distinguished reductive subalgebra of $V$; we 
denote it $\rr$.

Let $R$ contain a Virasoro element $L$ such that
\begin{equation}\label{eq.physpair}
[L_{\la}g]=(\d + \la)g \text{ for any } g \in \rr.
\end{equation}
Then $(R,L)$ is a \emph{simple physical Virasoro pair}.

By skew-symmetry, (\ref{eq.physpair}) is equivalent to
$[g_{\la}L]=\la g$, which implies that $[g_{(0)}L]=0$, i.e. that
$L$ is invariant with respect to the action of $\rr$ on $R$ (by
means of the $0$th product). Also, $R$ is completely reducible as
an $\rr$-module because the $0$th product commutes with the
$\Cset[\d]$-module structure.

\subsection{Classification results}
\begin{theorem}\lbb{veltsfslcs}
The following is a complete list of all simple physical Virasoro
pairs:
\begin{align*}
W_0 \quad :  \quad  &L=-1 \\
W_N,\; N \geq 1 \quad : \quad &L=-1 + (p_0 + p_1 \d)E  \\
S_{N,a},\; N \geq 2, \;  a \in \Cset \quad : \quad &L=-1 +
\frac{1}{N}(\d -a)E  \\
\wti{S}_N,\; N \geq 2, \; N \text{even} \quad : \quad &L=-(1
-\xi_1 \ldots \xi_N)
+ \frac{1}{N} \d E  \\
K_N,\;  N \geq 1, \; N \text{odd}  \quad : \quad &L=-1 \\
K_N, \; N>6, \; N \text{even} \quad : \quad &L=-1 + p_1 \d \nu \\
K_6 \quad:                   \quad &L=-1 + (p_1\d +p_3 \d^3)\nu   \\
K'_4 \quad:                    \quad &L=-1 + (p_0 + p_1\d)  \d \nu \\
CK_6 \quad :                    \quad &L=-1 + \alpha \d^3 \nu.
\end{align*}
where $E=\sum_{i=1}^N \xi_i \d_i$ is the Euler operator, 
$\nu=\xi_1\dots\xi_N$ is
the longest monomial in $\wedge(N)$, $p_0,p_1,p_3 \in \Cset$ and
$\alpha \in \Cset$, $\alpha^2=-1$.
\end{theorem}

\begin{proof}
Recall that $W_N=\Cset[\d]\otimes(W(N) \oplus \wedge(N))$. If
$N=0$, the equation $[L_{\la}L]=(\d + 2 \la)L$ implies that
$L=-1$. Suppose $N \geq 1$, then $\rr=\gl_N$.
First of all we remark that $L$ cannot be contained in the radical
of the even part of $W_N$, which we denote by $\Rad W_{N~\ov{0}}$.
We set
\begin{equation*}
S=\Cset[\d]\otimes(\oplus_{ k \geq 2, \; k \; \text{even}}W(N)^k),
\quad T=\Cset[\d]\otimes(\oplus_{k \geq 2, \; k \;
\text{even}}\wedge^k(N)).
\end{equation*}
Since the $\la$-bracket in $W_N$ is graded, $S\oplus T \subseteq
\Rad W_{N ~\ov{0}}$, hence any $L$ must have a component in
$\Cset[\d]\otimes (W(N)^0 \oplus \Cset 1)$:
\begin{equation*}
L=Q(\d)1 + P(\d)E + \sum_i P_i(\d) g_i + s +t,
\end{equation*}
where $\{ g_i \}$ is a basis of $\sl_N$, $s \in S$ and $t \in T$.
Furthermore, we can compare the coefficients of the degree $0$
basis elements in the equation $[L_{\la}L]=(\d + 2 \la)L$. In the
case of $1 \in \wedge(N)$ we have:
\begin{equation*}
-(\d + 2 \la)Q(-\la)Q(\d +\la)=(\d +2 \la)Q(\d),
\end{equation*}
hence $Q(\d)=Q$ is constant and either $Q=0$ or $Q=-1$. Comparing
the coefficients of $E$ we get:
\begin{equation*}
-(\d + \la)Q P(\d + \la) -\la Q P(-\la)=(\d + 2 \la)P(\d).
\end{equation*}
Comparing the terms with basis elements in $\sl_N$ we get
\begin{align*}
&(\d + \la)Q \sum_i P_i(\d + \la)g_i -\la Q \sum_i P_i(- \la)g_i
+ \sum_{i,j}P_i(-\la)P_j(\d + \la)[g_i, g_j] \\
&=(\d + 2 \la)\sum_i P_i(\d )g_i.
\end{align*}
Now, $Q=0$ implies $P=0$ and
\begin{equation*}
\sum_{i,j}P_i(-\la)P_j(\d + \la)[g_i, g_j]=(\d + 2 \la)\sum_i
P_i(\d )g_i.
\end{equation*}
If we examine the degree of $\d$ in the last equation, we see that
$P_i=0$ for any $i$. It follows that if $Q=0$, then $L=s+t \in
\Rad W_{N~\ov{0}}$, which is impossible. Hence $Q=-1$. If we
examine the equation for $E$, which is homogeneous, we conclude
that $P(\d)=p_0 + p_1 \d$ and $L=-1 + (p_0 + p_1 \d)E + \sum_i
P_i(\d) g_i + s +t$.

Next, we impose the condition that $[L_{\la}g]=(\d + \la)g$ for
any $g \in \gl_N$. Since $[-1_{\la}g]=(\d +\la)g$, we have
$[(\sum_i P_i(\d) g_i +s +t){}_{\la}g]=0$
for any $g\in \gl_N$.
Consequently, in degree $0$ we have $[(\sum_i P_i(\d)
g_i)_{\la}g]=0$ for any $g \in \gl_N$. In particular, it follows
that $\sum_i P_i(\d) g_i$ lies in the center of $\cur \sl_N$ which
is zero. On the other hand, in degree $\geq 2$ we have
$[(s+t)_{\la}g]=0$, hence $[g_{(0)}(s+t)]=0$. By looking at the
decomposition of the $\Cset[\d]$-basis of $W_N$ into irreducible
$\gl_N$-modules, we see that there are no invariants in degree
$\geq 2$. Consequently, $s+t=0$ and we conclude that
$L=-1 + (p_0 + p_1 \d)E$.

A similar argument, with $\sl_N$ replacing $\gl_N$, 
provides the solution for $S_{N,a}$ and $\wti{S}_N$.

The cases of $K_N$, $0 \leq N \leq 3$, $K'_4$ and $CK_6$ are dealt
with in a similar fashion. As for $K_N$, $N \geq 5$, the unique
$\so_N$-invariant is $\wedge^N(N)$. If $N$ is odd, this vector is
odd too, hence cannot appear in $L$. If $N$ is even, from
$[L_{\la}L]=(\d + 2 \la)L$ we obtain a homogeneous equation, which
admits a linear solution for any $N$ and a cubic solution for
$N=6$.

In the case of $\cur  \ss$, there are no non-zero Virasoro elements.
To see this, let $\{ g_i \}_{i \in I}$ be a linear basis of $\ss$.
Then $L=\sum_i P_i(\d)g_i$ and by substituting this expression into 
$[L_\la L]=(\d+2\la)L$ and comparing the degrees in $\d$, we see that 
$L=0$. 
\end{proof}

\begin{corollary}\lbb{evsedip}
The following is a complete list of all pairs $(R,L)$, where $R$ is a 
finite simple Lie conformal superalgebra and $L$ is a Virasoro element 
such that $R_{ \overline{0} }=L \sd \cur \g$, where
$\g$ is a finite-dimensional Lie algebra and $[L_\la a]=(\d+\la)a$ for all 
$a\in \g$:
\begin{align*}
W_1, \quad &-1+(p_0+p_1\d)E  \\
W_2, \quad &-1   + \frac{1}{2} \d E  \\
S_{ 2,a }, \quad &-1 + \frac{1}{2} (\d -a)E \\
\wti{S}_2, \quad &-(1 - \nu)  + \frac{1}{2} \d E \\
K_0, \quad &-1  \\
K_1,  \quad &-1  \\
K_3, \quad &-1  \\
K'_4, \quad &-1 +(p_0 + p_1 \d)\d \nu        \\
CK_6,  \quad &-1 +\alpha \d^3 \nu.
\end{align*}
\end{corollary}

\begin{proof}
We have to impose the condition $[L_{\la}g]=(\d + \la)g$ for any
$g\in \g$.

Let $w \in W(0,N)_k \subseteq W_N$ and $f \in \wedge^k(N)\subseteq
W_N$. We have
\begin{align*}
&[L_{\la}w]=(p_0 k + \d +(1-p_1 k)\la)w, \\
&[L_{\la}f]=(p_0 k + \d +(2-p_1 k)\la)f - \la (p_0 - p_1 \la)f E.
\end{align*}
If $N >2$, $W(0,N)_2 \neq 0$ and the first equation implies that
$1-2p_1=1$ i.e. $p_1=0$. On the other hand, $\wedge^2(N) \neq 0$
too and the second equation implies that $2-2p_1=1$ i.e.
$p_1=1/2$. The contradiction proves that $N \leq 2$. Since $W_0
\simeq K_0$ and $W_1 \simeq K_2$, we concentrate on $W_2=\Cset[\d]
\otimes(W(2) \oplus \wedge(2))$. If we apply the second equation
to $ \xi_1 \xi_2 $, we get  $ p_0=0 $ and $ p_1= 1/2 $.

If $R=S_{ N,a }$ or $R=\wti{S}_N$, then $R=\Cset[\d] \otimes W(N)$
as $\sl_N$-modules. If $N >2$, $W(N)^2 \neq 0$ and the formulas
that appear in the proof of \cite[Proposition 4.16]{FK} show that
the conformal weight of any vector in $W(N)^2$ is either $(N-2)/N$ or $(2N 
-2)/N \neq 1$. Therefore $N=2$ is the only possibility.

If $R=K_N$, $N \geq 5$, the conformal weight of $\xi_1 \xi_2 \xi_3
\xi_4$ with respect to any simple physical Virasoro pair $(K_N,L)$
is $0$.

As for $K'_4$, we can check by direct computation that $K'_{4 \;
\overline{0}}=L \sd Cur \; \cso_4$. Similarly in the case of
$K_N$, $0 \leq N \leq 3$.

Finally, it was proved in \cite{CK2} that $CK_{6 \;\overline{0} }
=L \sd Cur \; \so_6$.
\end{proof}

\begin{remark}
For $W_2$ we have $\g=\sl_2 \oplus \mathfrak{b}_2$; here
$\mathfrak{b}_2$ denotes the two-dimensional solvable Lie algebra.
In all other cases, $\g=\rr$.
\end{remark}

\subsection{Physical conformal superalgebras}
A \emph{physical} Lie conformal superalgebra is a simple physical
Virasoro pair $(R,L)$ such that the following conditions hold:
\begin{enumerate}
\item
$R_{\ov{0}} =L \sd \cur  \g$, where $\g$ is a finite-dimensional
Lie algebra and $\rr\subseteq \g$;
\item
$L_{(0)}a = \d a$ for any $a \in R$.
\end{enumerate}

We thus obtain a generalization of the main result of \cite{K4} (see also 
\cite{Y}).

\begin{proposition}\label{physlist}
The following is a complete list of physical Lie conformal
superalgebras:
\begin{align*}
W_1, \quad &-1+p_1\d E \\
W_2, \quad &-1  + \frac{1}{2} \d E  \\
S_{2,0}, \quad &-1 + \frac{1}{2} \d E \\
K_0, \quad &-1  \\
K_1,  \quad &-1  \\
K_3, \quad &-1  \\
K'_4, \quad &-1 +(p_0 + p_1 \d)\d \nu        \\
CK_6,  \quad &-1 +\alpha \d^3 \nu.
\end{align*}
\end{proposition}

\begin{proof}
Condition (2) excludes $\wti{S}_2$ from the list in Corollary
\ref{evsedip}, implies that $a=0$ for $S_{2,a}$ and forces $p_0=0$
for $(W_1,L)$.
\end{proof}

The formal distribution Lie superalgebras $(\Lie R, R)$ corresponding to 
the Lie conformal superalgebras on the list of Proposition~\ref{physlist} 
are well known to physicists (except for $CK_6$).  All of them, except for 
$CK_6$, have non-trivial central extensions.  In the cases 
$R=K_0=W_0,K_1,K_2=W_1,K_3,S_{2,0},$ and $K_4'$, these central extensions 
are 
called the Virasoro, Neveu--Schwarz, $N=2$, $N=3$, $N=4$, and big $N=4$ 
superconformal algebras.  With the exception of the last one, they appear 
on the list of \cite{RS}.  The last one appeared independently in 
\cite{KL}, \cite{S}, and \cite{STP}.


\section{Representations of Solvable Lie Conformal Superalgebras}\label{sec.lieth}

In this section we discuss several results towards a possible
analog of Lie Theorem for finite Lie conformal superalgebras.

\subsection{Preliminaries} First we restate the analogs of Lie theorems in 
the, respectively, non-con\-formal and conformal non-super cases:

\begin{theorem}\label{nonconflie}\cite{K3} Let $L$ be a finite-dimensional
solvable Lie superalgebra and $V$ its irreducible representation.  Then
either $\dim V_{\ov 0}=\dim V_{\ov 1}$ and $\dim V=2^s$, where $s\leq \dim
L_{\ov 1}$ or $\dim V=1$.
\end{theorem}

\begin{theorem}\label{nonsuplie}\cite{DK} Let $R$ be a finite solvable Lie
conformal algebra and $M$ an $R$-module.  Then there exists $v\in
M$ such that $x_\la v=\varphi(a)v$, $\varphi: R\to\CC[\la]$.  In
particular, every non-trivial simple $R$-module $M$ is free
and $\rk M=1$.
\end{theorem}

\begin{remark} Let $R$ be a Lie conformal algebra with a module $M$.
Then $(\Tor R)_\lambda M=0$, $\Tor M$ is a submodule of $M$, and $R_\la
\Tor M=0$ \cite{DK}.  Since we are interested in simple modules only,
below we will always  assume that both the conformal algebra and its
modules are free over $\CC[\d]$.
\end{remark}

We will also need the following

\begin{lemma}\cite{DK}\label{rankdown} For a solvable Lie conformal
superalgebra $R$, $\rk
R>\rk R'$, where $R'=\<R_\la R\>$ is the derived algebra of $R$.
\end{lemma}

Another useful result is the analog of Lemma~\ref{milcs}:

\begin{lemma}\label{milcsmod} A finite module $M$ over a finite Lie
conformal
superalgebra $R$ either contains a simple $R$-submodule or a
submodule $N$ such that $R_\la N=0$.
\end{lemma}

\subsection{Rank $1$ modules}
Let $R$ be a finite solvable Lie conformal superalgebra and $M$ an
$R$-module, $\rk M=1$.  In particular, $M=\CC[\d]v$ for some $v\in M$ and
the representation of $R$ is given by the function $\ell: R\to\CC[\la,\d]$
such that $x_\la v=\ell(x)v$ for all $x\in R$.  Theorem~\ref{nonsuplie}
immediately implies that either $\ell=0$ or $M$ is a simple $R_{\ov
0}$-module and that $\ell$ as a function on $R_{\ov 0}$ depends only on
$\la$.  A direct
computation shows that $\ell(R_{\ov 0}')=0$. Also, since the action of
$R_{\ov 1}$ must change parity, we conclude that $\ell(R_{\ov 1})=0$.

Conversely, if $\ell$ satisfies all the above conditions and $\ell(\d
x)=-\la \ell(x)$ for all $x$, $\ell$ determines an $R$-module of rank one.

Let $\mathcal{L}$ be the set of functions $\ell:R\to\CC[\la]$ such that
$\ell(R_{\ov 0}')=\ell(R_{\ov 1})=0$ and $\ell(\d x)=-\la \ell(x)$ for all
$x\in R$.  Let $\mathcal{L}_0$ be a subspace of $\mathcal{L}$ consisting
of $\ell$ for which $\ell(R')=0$.

\begin{remark}  Let $H$ be a subalgebra of $R$ such that $R=H\oplus
\CC[\d]x$, $x\in R$, $R_{\ov 0}\subset H$.   The $\la$-bracket in $R$
induces an $H$-action on $\CC[\d]x$.  The corresponding $\ell$ lies in
$\mathcal{L}_0$.
\end{remark}

Let $M$ be an $R$-module and $\ell\in\mathcal{L}_0$.  Let $V$ be a
$\CC[\d]$-generating subspace of $M$, i.e. $M=\CC[\d]\otimes V$.  We
introduce another $R$-module structure on $M$ by putting $x_\la v=x^M_\la
v+\ell(x)v$ for $x\in R$, $v\in V$ and extending the $R$-action to all of
$M$ in the standard way.  In this way we obtain a module $M'$ and call $M$
and $M'$ \emph{$\mathcal{L}_0$-equivalent}.

\subsection{Induced modules} In the proof of
Theorem~\ref{nonconflie}, the powers of $2$ appear because, in
principle, a simple module $V$ over a Lie superalgebra $L$
can be induced from a smaller module $W$ over a subalgebra $H$ of
$L$. In particular, when $L=H\oplus\CC g$ and $p(g)=\ov 1$,
$V=W+gW$, hence the dimension gets doubled.

In the conformal case, non-trivial induction is still possible but
the rank can grow arbitrarily as the following lemma demonstrates.

\begin{lemma}\label{ind1} Let $R$ be a Lie conformal superalgebra, and let
$H$ be a (homogeneous) subalgebra of $R$ such that $R'\subset H$ and
$R=H\oplus\CC[\d]x$ for some homogeneous $x\in R$.  Let $M$ be a finite
$R$-module and $N$ an $H$-submodule of $M$ that generates $M$. Then
\begin{equation}\label{eq.ind}
M=\sum_{i=1}^l \CC[\d]\left(\prod_{j=i}^l x\encirc{n_j}\right)N
\text{ for some } n_1,\dots,n_l.
\end{equation}
Moreover, as $H$-modules $M$ and $N$ have the same simple quotients.
If $N$ is simple, then $\rk M$ is proportional to $\rk N$.
\end{lemma}

\begin{proof}
Clearly, $\<x_\la N\>+N$ is an $H$-submodule of $M$ and, unless $M=N$, it
contains $N$ as a proper submodule.

Hence $U=(\<x_\la N\>+N)/N$ is an $H$-module.  Moreover, since
$x\encirc{k}(h\encirc{m}w)=\pm h\encirc{m}(x\encirc{k}w) \mod N$
for $k,m\in\ZZ_+$ and $w\in W$, we obtain submodules
$U_k=(\CC[\d](x\encirc{k}N)+N)/N$ of $U$.  The Leibniz rule
implies that $U_k\subset U_{k+1}$.  Obviously, $U=\sum U_k$ and it
immediately follows that $U=U_n$ for some $n$.  Hence, $N+\<x_\la
N\>=N+\CC[\d](x\encirc{n}N)$.  Now we can consider the latter
module instead of $N$.  Then (\ref{eq.ind}) follows by induction.

We have actually obtained a filtration of $M$: $N\subseteq
N+\CC[\d](x\encirc{0}N)\subseteq \ldots \subseteq M$.  Denote the
$k$-th element of the fil\-tra\-tion by $M_k$.  Then a direct
compu\-ta\-tion shows that there exists a natural surjective map $N\to
M_k/M_{k-1}$ for every $k$ (e.g. $N\to U_k/U_{k-1}$ is defined as
$v\mapsto x\encirc{k}v$).  This proves the rest of the lemma.
\end{proof}

\begin{remark} By the Jacobi identity, one can always have $n_j\leq 
n_{j+1}$
in (\ref{eq.ind}).  When $x$ is odd the inequality is strict.  The
analogous lemma for ordinary  algebras would then imply that $l=1$
whenever $x$ is odd but this is not true for conformal algebras.
\end{remark}

Let $H$ be a subalgebra of $R$, $M$ a simple $R$-module and $N$ a
simple $H$-submodule of $R$ such that $N$ generates $M$ over $R$.
Then we say that $M$ is \emph{induced} from $N$.  This does not mean that
we can define an induction functor from $H$-modules to $R$-modules; the
word ``induced'' is used here only as a shorthand.

Thus Lemma~\ref{ind1} provides a description of induced modules when
$\rk H=\rk R-1$ and $R'\subset H$.

\begin{proposition}  Let $M$ be a finite simple module over a finite
solvable Lie conformal superalgebra $R$.  Then all simple factors of
$M$ considered as an $R_{\ov 0}$-module are of rank $1$.  The
corresponding elements of $\mathcal{L}$ extended by zero to $R_{\ov 1}$
lie in a single coset $\ell_M\in\mathcal{L}/\mathcal{L}_0$.
\end{proposition}

\begin{proof} We use induction on $\rk R$.

Let $H$ be a subalgebra of $R$ such that $R=H\oplus\CC[\d]x$, $R'\subset
H$, and $x$ is homogeneous.  (Such $H$ always exists by
Lemma~\ref{rankdown}.) Two cases can occur: $M$ contains a simple
$H$-module $N$ or there exists a submodule $N\subset M$ such that $H_\la
N=0$.

Consider the first case, i.e. let $M$ be induced from $N$.  Then by
Lemma~\ref{ind1} all $H$-simple factors of $M$ are isomorphic to $N$.
Thus, if $x$ is odd, $M$ and $N$ viewed as $R_{\ov 0}$-modules have the
same simple factors and we are done by induction.   However, if $x$
is even, $M$ and $N$ have the same simple factors only as $H_{\ov
0}$-modules.  However, let then $\ell_1,\ell_2$ be the two forms
corresponding to the  simple factors and let $\ell=\ell_1-\ell_2$.  By
induction $\ell(R_{\ov 1}')=0$, hence we are done.

Consider now the second case, i.e. let $M$ be induced from $N$
such that $H_\la N=0$.  Then all simple factors of $M$ are
killed by $H$. Hence, if $x$ is odd, $M$ has no $R_{\ov
0}$-simple factors. If $x$ is even and we get a simple
factor that corresponds to $\ell\in\mathcal{L}$, we get
$\ell(R_{\ov 1}')=0$ and thus $\ell\in\mathcal{L}_0$.
\end{proof}



\end{document}